%% file: lls-non-ct-trop.tex
\documentclass{amsart}
\include{headers}

\begin{document}
\title{Limit linear series and the Amini-Baker construction}
\author{Brian Osserman}
\begin{abstract} We draw comparisons between the author's recent 
construction of limit linear series for curves not of compact type and
the Amini-Baker theory of limit linear series on metrized complexes,
as well as the related theories of divisors on discrete graphs and
on metric graphs. From these we conclude that the author's theory
(like the others) satisfies the Riemann and Clifford
inequalities. Motivated by our comparisons, we also develop negative 
results on Brill-Noether generality for certain families of 
metric graphs. Companion work of He develops our comparisons further
and uses them to prove new results on smoothability of Amini-Baker
limit linear series and of divisors on metric graphs.
\end{abstract}

\thanks{The author was partially supported by 
Simons Foundation grant \#279151 during the preparation of this work.}
\maketitle


\section{Introduction}

In \cite{os25}, the author introduced a theory of limit linear series
for nodal curves not of compact type. This was further studied in
\cite{os23}, yielding some suggestive connections to the tropical
proof of the Brill-Noether theorem by Cools, Draisma, Robeva and 
Payne \cite{c-d-p-r}. Separately, Amini and Baker \cite{a-b1} had
introduced an alternate notion of limit linear series for curves not
of compact type, more closely connected to the recently developed 
theories of divisors on discrete graphs \cite{b-n1} and on abstract 
tropical curves \cite{g-k2}, \cite{m-z1}.
The purpose of this paper is to examine various aspects of the 
connections between these theories. 

We begin by verifying some basic compatibilities between our notions
of multidegrees on dual graphs and the theory of divisors on (discrete)
graphs. As an application, we can leverage the existence of $v$-reduced
divisors on graphs to prove in Theorem \ref{thm:riemann} our own version 
of the ``Riemann's theorem'' proved both 
by Eisenbud and Harris (Theorem 4.1 of \cite{e-h1}) and Amini and Baker
(Remark 5.8 of \cite{a-b1}) for their respective limit linear series 
theory. This states in particular that when the degree $d$ is greater than 
$2g-2$, limit linear series only exist when $r \leq d-g$.

We next show in Theorem \ref{thm:forget} that there is a forgetful map 
between our limit linear series and Amini-Baker limit linear series.
Because Amini and Baker proved Theorem \ref{thm:riemann} as well as
Clifford's inequality for their limit linear series, we immediately
obtain in Corollary \ref{cor:riemann-clifford} both a second proof of 
Theorem \ref{thm:riemann}, and a proof that our limit linear series
satisfy Clifford's inequality.

There is a class of curves -- those of `pseudocompact type'\footnote{These
are the curves for which, if you squint hard enough, the dual graph is a
tree; see Definition \ref{def:pseudocompact} below.} -- for which
we have an alternative definition of limit linear series in terms of
generalizing vanishing conditions and gluing conditions, which is equivalent 
to (but formulated quite differently from) our more general definition. We 
show in Proposition \ref{prop:forget-pseudocompact} that
for these curves, the construction of our forgetful map will yield an 
Amini-Baker limit linear series even if our gluing condition is not 
satisfied. In a companion paper, Xiang He \cite{he2} shows that conversely, 
an Amini-Baker limit linear series on a curve of pseudocompact type satisfies 
our generalized vanishing condition, so that on such curves, Amini-Baker
limit linear series are equivalent to tuples of linear series which 
satisfy our generalized vanishing condition. He also examines cases
in which the forgetful map is and is not surjective, proving new 
results on smoothability (and non-smoothability) of Amini-Baker limit
linear series, with some applications also to smoothability of tropical
linear series.

Finally, in \S \ref{sec:examples} we give some examples 
highlighting the differences between our theory, the Amini-Baker theory,
and the theory of divisors on graphs. These include examples of curves 
which have maximal gonality in our sense but are hyperelliptic in the
Amini-Baker sense, and examples of curves which are not hyperelliptic
either in our sense or the Amini-Baker sense, but which carry a 
divisor of degree $2$ and rank $1$ on the underlying metric graph.
Motivated by these examples, and in a similar vein to the recent work of
Kailasa, Kuperberg and Wawrykow \cite{k-k-w1}, we then develop more 
systematic negative results on Brill-Noether generality for graphs,
showing that if a graph has a point disconnecting it into three or more
components, or is obtained by attaching two graphs to one another by a
collection of four or more edges with the same endpoints, then the graph
is not Brill-Noether general. We also 
discuss the overall philosophy of when one could reasonably expect a graph 
to be Brill-Noether general.

We also include an appendix with some background results on the Amini-Baker
theory, especially relating to restricted rank.

Because both our limit linear series of \cite{os25} and the Amini-Baker
limit linear series use $\Gamma$ and $G$ in different ways as basic 
notation, we are not able to follow both at once. We have decided to follow 
the notation of \cite{os25} as is, while using $\GamAB$ and $\GAB$ for
the Amini-Baker usage of $\Gamma$ and $G$, respectively.

\subsection*{Acknowledgements} 
I would like to thank Omid Amini, Matt Baker, Vivian Kuperberg and Sam Payne 
for several helpful conversations, particularly regarding the material in 
the final section and the appendix. I would also like to thank Xiang He for 
many helpful comments and conversations during the course of preparing this 
work.

\section{Background on limit linear series}\label{sec:defs}

In this section, we recall the general definition of limit linear series 
introduced in \cite{os25}. We also develop some new definitions and
results at the end of the section.

We begin with some definitions of a combinatorial nature. For us,
a \textbf{multidegree} on a graph is simply an integer vertex weighting
(i.e., exactly what is called a divisor in Brill-Noether theory for
graphs).
However, we will be interested in a notion of ``admissible'' multidegrees,
in which we distinguish between ``original'' vertices and vertices 
introduced by subdivision of edges. While this introduces some extra
complications into the notation, it is ultimately an important simplifying
tool, especially in the context of curves of pseudocompact type
(discussed in \S \ref{sec:forget-pct}). For instance, in the special
case of curves of two components, it means that our limit linear series
will always only involve a pair of linear series on components, even
when it may have been necessary to introduce new rational components 
to extend the underlying line bundle.

In the below,
$\Gamma$ will be obtained by choosing a directed structure on the dual
graph of a projective nodal curve.

\begin{defn} A \textbf{chain structure} on a graph $\Gamma$ is a function
$\bn:E(\Gamma) \to \ZZ_{> 0}$. A chain structure is \textbf{trivial}
if $\bn(e)=1$ for all $e \in E(\Gamma)$. Given $\Gamma$ and $\bn$, let
$\widetilde{\Gamma}$ be the graph obtained from $\Gamma$
by subdividing each edge $e$ into $\bn(e)$ edges. 
\end{defn}

Thus, we have a natural inclusion $V(\Gamma) \subseteq V(\widetilde{\Gamma})$,
and each vertex $v$ of $\widetilde{\Gamma}$ not in $V(\Gamma)$ is naturally
associated to an edge $e \in E(\Gamma)$; we will refer to $v$ as a 
``new vertex lying over $e$.'' 
The chain structure will determine the length of the chain of rational
curves inserted at a given node, so that $\widetilde{\Gamma}$ will be the
dual graph of the resulting curve. Note that the trivial case (in which no 
rational curves are inserted) corresponds to $\bn(e)=1$. 

As may be
suggested by the construction of $\widetilde{\Gamma}$, the analogue of
our chain structure in the tropical setting is the edge lengths inducing
a metric graph structure on a given dual graph.

\begin{defn} If $\Gamma$ is a directed graph, 
for each pair of an edge $e$ and adjacent vertex $v$ of $\Gamma$,
let $\sigma(e,v)=1$ if $e$ has tail $v$, and $-1$ if $e$ has head $v$.
Given also $\bn$ a chain structure on $\Gamma$, an
\textbf{admissible multidegree} $w$ of total
degree $d$ on $(\Gamma,\bn)$ consists of a function
$w_{\Gamma}: V(\Gamma) \to \ZZ$ together with a tuple
$(\mu(e))_{e \in E(\Gamma)}$, where each $\mu(e) \in \ZZ/\bn(e)\ZZ$,
such that
$$d = \#\{e \in E(\Gamma): \mu(e) \neq 0\}
+ \sum_{v \in V(\Gamma)} w_{\Gamma}(v).$$

The multidegree $\widetilde{w}$ on $\widetilde{\Gamma}$ induced by $w$
is defined by $\widetilde{w}(v)=w_{\Gamma}(v)$ for
all $v \in V(\Gamma)$, by $\widetilde{w}(v)=1$ if $\mu(e)\neq 0$ and
$v$ is the $\mu(e)$th new vertex lying over $e$, and by
$\widetilde{w}(v)=0$ otherwise. Here we order the new vertices over a
given $e \in E(\Gamma)$ using the direction of $e$.
\end{defn}

The idea behind admissible multidegrees is that if we have a line bundle
on the generic fiber of a one-parameter smoothing of a nodal curve, in 
order to extend it to the special fiber, it suffices to consider 
multidegrees which have degree $0$ or $1$ on each rational curve inserted 
at the node, with degree $1$ occurring at most once in each chain. Thus, 
$\mu(e)$ determines where on the chain (if anywhere) positive degree occurs.

\begin{defn} Given a chain structure $\bn$ on $\Gamma$,
let $w$ be an admissible multidegree. Given also
$v \in V(\Gamma)$, the \textbf{twist} of $w$ at $v$ is obtained as follows:
for each $e$ adjacent to $v$, increase $\mu(e)$ by $\sigma(e,v)$. Now,
decrease $w_{\Gamma}(v)$ by the number of $e$ for which $\mu(e)$ had
been equal to $0$, and for each $e$, if the new $\mu(e)$ is zero, increase
$w_{\Gamma}(v')$ by $1$, where $v'$ is the other vertex adjacent to $v$.
The \textbf{negative twist} of $w$ at $v$ is the admissible multidegree
$w'$ such that the twist of $w'$ at $v$ is equal to $w$.
\end{defn}

Twists will be the change in multidegrees accomplished by twisting by
certain natural line bundles; see Notation \ref{not:twist} below.

\begin{ex}\label{ex:twists-trivial} In the case of trivial chain structure,
a twist at $v$ simply reduces $w_{\Gamma}(v)$ by the valence of
$v$ while increasing $w_{\Gamma}(v')$ by the number of edges connecting
$v'$ to $v$, for each $v' \neq v$. This is the same as the chip firing
considered by Baker and Norine in \cite{b-n1}.
\end{ex}

\begin{ex}
\begin{figure}
\centering
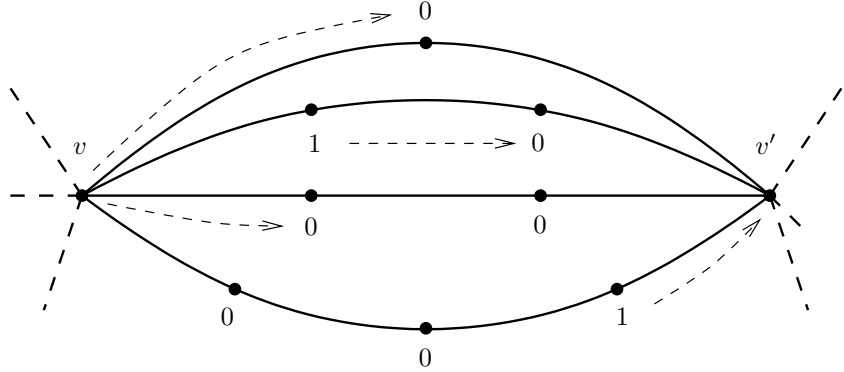
\caption{An admissible multidegree and the effect of twisting at $v$.}
\label{fig:twist}
\end{figure}
The effect of a twist at $v$ is shown in Figure \ref{fig:twist}, which
shows adjacent vertices $v,v'$ of $\Gamma$, connected by four edges
in $\Gamma$, which are subdivided in $\widetilde{\Gamma}$ according to
the chain structure. The effect of twisting at $v$ is then that for each
chain between $v,v'$ in the figure, the $1$ is moved one vertex to the
right. When there is no $1$ in the chain, a $1$ is placed on the first
vertex, and degree at $v$ is reduced by $1$. When there is a $1$ on the
last vertex, it is removed and the degree at $v'$ is increased by $1$.
Thus, for the example in the figure, twisting at $v$ will decrease the
degree at $v$ by $2$, and increase the degree at $v'$ by $1$.
\end{ex}

\begin{rem}\label{rem:twists} Induced multidegrees on $\widetilde{\Gamma}$
are compatible with twists as follows:
twisting $w$ at $v \in V(\Gamma)$ is
the same as twisting $\widetilde{w}$ by $v$, and then also by all
new vertices between $v$ and the $\sigma(e,v)\mu(e)$th new vertex lying over
$e$, for each $e \in E(\Gamma)$ adjacent to $v$.
\end{rem}

\begin{defn} An admissible multidegree $w$ is \textbf{concentrated} at
a vertex $v \in V(\Gamma)$ if there is an ordering on $V(\Gamma)$ starting
with $v$, and such that for each subsequent vertex $v'$, we have that $w$
becomes negative in index $v'$ after taking the composition of the negative
twists at all the previous vertices.
\end{defn}

We will also refer to a multidegree on a graph without chain structure
as being concentrated at $v$ if it is concentrated at $v$ when considered
as an admissible multidegree with respect to the trivial chain structure.

While the definition of concentrated might be a bit opaque, it is based
on a very simple geometric concept: if we have a line bundle
$\sL$ of multidegree $w$ on a curve with dual graph $\Gamma$, the definition
of concentrated implies that if a section of $\sL$ vanishes on the 
component $Z_v$ corresponding to $v$, then we can iteratively traverse
the other components to conclude it must vanish everywhere. We will show
in Corollary \ref{cor:vred-conc} that concentratedness is essentially
equivalent to the notion of $v$-reducedness for divisors on graphs, except
that nonnegativity away from $v$ is imposed for the latter. See
Remark \ref{rem:v-reduced-role} for discussion of why we adopt the more
general condition.


The following directed graph keeps track of all the multidegrees we will
want to consider starting from any one admissible multidegree.

\begin{notn} Let $G(w_0)$ be the directed graph with vertex set
$$V(G(w_0)) \subseteq \ZZ^{V(\Gamma)} \times \prod_{e \in E(\Gamma)}
\ZZ/\bn(e)\ZZ$$
consisting of all admissible multidegrees obtained from $w_0$ by sequences
of twists, and with an edge from $w$ to $w'$ if $w'$ is obtained from $w$
by twisting at some vertex $v$ of $\Gamma$.

Given $w \in V(G(w_0))$ and $v_1,\dots,v_m \in V(\Gamma)$ (not necessarily
distinct), let $P(w,v_1,\dots,v_m)$ denote the path in $V(G(w_0))$ obtained
by starting at $w$, and twisting successively at each $v_i$.
\end{notn}

By the invertibility of
twists, $G(w_0)=G(w)$ if and only if $w \in G(w_0)$. While our directed
structure on $\Gamma$ is just a convenience, the directedness of $G(w_0)$
is crucial. 
Also, note that the endpoint of $P(w,v_1,\dots,v_m)$ is independent of the
ordering of the $v_i$. In fact, we have the following (see Proposition 2.12
of \cite{os25}, although this is also standard in the chip-firing 
literature; see for instance Lemma 2.2 of \cite{h-l-m-p-p-w1}):

\begin{prop}\label{prop:paths-unique} If $P(w,v_1,\dots,v_m)$ is a minimal
path in $G(w_0)$ from $w$ to some $w'$, then $m$ and the $v_i$ are
uniquely determined up to reordering.

More generally, paths $P(w,v_1,\dots,v_m)$ and $P(w,v'_1,\dots,v'_{m'})$
have the same endpoint if and only if the multisets of the $v_i$ and the
$v'_i$ differ by a multiple of $V(\Gamma)$.
\end{prop}

We now move on to definitions which involve geometry more directly.

\begin{sit}\label{sit:basic}
Let $X_0$ be a projective nodal curve, with dual graph
$\Gamma$, and choose an orientation on $\Gamma$. For $v \in V(\Gamma)$,
let $Z_v$ be the corresponding irreducible component of $X_0$, and $Z_v^c$
the closure of the complement of $Z_v$ in $X_0$.
\end{sit}

A preliminary definition is the following.

\begin{defn}\label{def:enriched} If $X_0$ is a nodal curve with dual
graph $\Gamma$, an
\textbf{enriched structure} on $X_0$ consists of the data, for each
$v \in V(\Gamma)$ of a line bundle $\sO_v$ on $X_0$,
satisfying the following conditions:
\begin{Ilist}
\itm for any $v \in V(\Gamma)$, we have
$$\sO_v|_{Z_v} \cong \sO_{Z_v}(-(Z_v^c \cap Z_v)),\text{ and }
\sO_v|_{Z_v^c} \cong \sO_{Z_v^c}(Z_v^c \cap Z_v);$$
\itm we have
$\bigotimes_{v \in V(\Gamma)} \sO_v \cong \sO_{X_0}$.
\end{Ilist}
\end{defn}

Enriched structures are induced by regular one-parameter smoothings
$\pi:X \to B$, by setting $\sO_v=\sO_X(Z_v)|_{X_0}$. They are necessary
data for our definition of limit linear series, but because they do not 
occur in the Amini-Baker or tropical settings, they will play a relatively
minor role in the present paper.

We now explicitly introduce the chains of rational curves induced by a
chain structure on $X_0$.

\begin{defn}
Given $X_0$ and a chain structure $\bn$ on $\Gamma$, let $\widetilde{X}_0$ 
denote the nodal curve obtained from $X_0$ by, for each $e \in E(\Gamma)$,
inserting a chain of $\bn(e)-1$ projective lines at the corresponding
node. We refer to the new components of $\widetilde{X}_0$ as the
\textbf{exceptional components}.
\end{defn}

Thus, $\widetilde{\Gamma}$ is the dual graph of $\widetilde{X}_0$,
and an admissible multidegree on $\Gamma$ induces a usual multidegree on
$\widetilde{X}_0$.

From now on, we will assume we have fixed an enriched structure together
with suitable global sections, as follows.

\begin{sit}\label{sit:enriched} In Situation \ref{sit:basic}, suppose
we have also a chain structure $\bn$ on $\Gamma$, and an enriched
structure $(\sO_v)_v$ on the resulting $\widetilde{X}_0$, and for
each $v \in V(\widetilde{\Gamma})$, fix
$s_v \in \Gamma(\widetilde{X}_0,\sO_v)$ vanishing precisely on $Z_v$.
\end{sit}

The sections $s_v$ will be convenient in describing maps between different
twists of line bundles; they will not be unique even for curves of compact
type, but they do
not ultimately affect our definition of limit linear series. See
Remark 2.22 of \cite{os25}.

We next describe how, given an enriched structure on $\widetilde{X}_0$, and
a line bundle $\sL$ of multidegree $w_0$, we get a collection of line bundles
indexed by $V(G(w_0))$, with morphisms between them indexed by $E(G(w_0))$.

\begin{notn}\label{not:twist} In Situation \ref{sit:enriched} assume we
are given also an admissible multidegree $w_0$ on $(\Gamma,\bn)$. Then
for any edge $\e \in E(G(w_0))$, starting at
$w=(w_{\Gamma},(\mu(e))_{e \in E(\Gamma)})$ and determined by twisting at
$v \in V(\Gamma)$, we have the corresponding twisting line bundle
$\sO_{\e}$ on $\widetilde{X}_0$ defined as
$$\sO_{\e}=\sO_v \otimes \bigotimes_{e \in E(\Gamma)}
\bigotimes_{i=1}^{\sigma(e,v)\mu(e)} \sO_{v_{e,i}},$$
where the first product is over edges $e$ adjacent to $v$, and for any
such pair, $v_{e,i}$ denotes the $i$th rational curve in
$\widetilde{X}_0$ from $Z_v$ on the chain corresponding to $e$.

In addition, we have the section $s_{\e}$ of $\sO_{\e}$ obtained from
the tensor product of the relevant sections $s_v$ and $s_{v_{e,i}}$.

Similarly, given $w,w' \in V(G(w_0))$, let $P=(\e_1,\dots,\e_m)$ be
a minimal path from $w$ to $w'$ in $G(w_0)$, and set
$$\sO_{w,w'}=\bigotimes_{i=1}^m \sO_{\e_i}.$$
\end{notn}

In Notation \ref{not:twist}, we take the
representative of $\sigma(e,v)\mu(e)$ between $0$ and $\bn(e)-1$, and
if $\mu(e)=0$, the product over $i$ is empty for the given $e$.
Note that it follows from Proposition \ref{prop:paths-unique} that
the constructions of Notation \ref{not:twist} are independent of choices
of (minimal) paths. The reason for the notation $\sO_{w,w'}$ is that,
as one can easily verify, tensoring by $\sO_{w,w'}$ take a line bundle of
multidegree $w$ to one of multidegree $w'$.

\begin{notn}\label{not:more-twist}
In Situation \ref{sit:enriched}, suppose $\sL$ is a line bundle on
$\widetilde{X}_0$ of multidegree $w_0$. Then for any $w \in V(G(w_0))$, set
$$\sL_w := \sL \otimes \sO_{w_0,w}.$$

Given also $w_v \in V(G(w_0))$ concentrated at $v$, set
$$\sL^v:=\sL_{w_v}|_{Z_v}.$$

Given an edge $\e$ from $w$ to $w'$ in $G(w_0)$, corresponding to
twisting at $v$, then either
$\sL_{w'}=\sL_w \otimes \sO_{\e}$, or
$\sL_w = \sL_{w'} \otimes \sO_{w',w}$. In the former case, we get a morphism
$\sL_w \to \sL_{w'}$ induced by $s_{\e}$. In the latter case, we observe
that $\sO_{w',w} \otimes \sO_{\e} \cong \sO_{\widetilde{X}_0}$, and
fixing such an isomorphism and again using $s_{\e}$ gives an induced
morphism
$$\sL_w \to \sL_w \otimes \sO_{\e} =
\sL_{w'}\otimes \sO_{w',w} \otimes \sO_{\e} \cong \sL_{w'}.$$
In either case, pushing forward gives an induced morphism
$$f_{\e}:\Gamma(\widetilde{X}_0, \sL_w) \to
\Gamma(\widetilde{X}_0, \sL_{w'}).$$

Finally, if $P=(\e_1,\dots,\e_m)$ is any path in $G(w_0)$, set
$$f_P:=f_{\e_m} \circ \cdots \circ f_{\e_1}.$$
If $P$ is a minimal path from $w$ to $w'$, write
$$f_{w,w'}:=f_P.$$
\end{notn}

We have the following simple consequence of Proposition
\ref{prop:paths-unique}:

\begin{cor}\label{cor:maps-defined} For any $w,w' \in V(G(w_0))$, the
morphism $f_{w,w'}$ is independent of the choice of minimal path.
\end{cor}

We can now give the definition of a limit linear series. 

\begin{defn}\label{def:lls}
Let $X_0$ be a projective nodal curve, $\bn$ a chain structure,
$w_0$ an admissible multidegree of total degree $d$ on $(X_0,\bn)$, and
$(\sO_v)_{v \in V(\Gamma)}$ an enriched structure on $\widetilde{X}_0$.
Choose also a tuple $(w_v)_{v \in V(\Gamma)}$ of vertices of $G(w_0)$, with
each $w_v$ concentrated at $v$, and sections $(s_v)_v$ as in Situation
\ref{sit:enriched}.  Then a \textbf{limit linear series}
on $(X_0,\bn)$ consists of a line bundle $\sL$ of multidegree $w_0$ on
$\widetilde{X}_0$, together with $(r+1)$-dimensional subspaces $V^v$ of 
$\Gamma(Z_v,\sL^v)$ for
each $v \in V(\Gamma)$, satisfying the condition that for all
$w \in V(G(w_0))$, the natural morphism
\begin{equation}\label{eq:gluing-map-1}
\Gamma(\widetilde{X}_0, \sL_w)
\to \bigoplus_{v \in V(\Gamma)} \Gamma(Z_v, \sL^{v})/V^v
\end{equation}
has kernel of dimension at least $r+1$, where \eqref{eq:gluing-map-1} is
obtained as the composition
\begin{multline*}
\Gamma(\widetilde{X}_0, \sL_w) \overset{\oplus f_{w,w_v}}{\to}
\bigoplus_{v \in V(\Gamma)} \Gamma(\widetilde{X}_0, \sL_{w_v})
\\ \to \bigoplus_{v \in V(\Gamma)} \Gamma(Z_v, \sL^{v})
\to \bigoplus_{v \in V(\Gamma)} \Gamma(Z_v, \sL^{v})/V^v.
\end{multline*}
\end{defn}

Although the choices of concentrated multidegrees are necessary to
even define the data underlying a limit linear series, any two choices
give canonically isomorphic moduli spaces; see Proposition 3.5 of \cite{os25}.

We now develop some new material which will be important in comparing to
the Amini-Baker theory, and should in any case be of independent interest.
It relates to the following bounded subgraph of multidegrees, generalizing 
the construction given in the compact type case in Definition 3.4.9 of 
\cite{os20}.

\begin{notn}\label{notn:barg}
In the situation of Definition \ref{def:lls},
let $\bar{G}(w_0)$ be the subgraph of $G(w_0)$ obtained
by restricting to $w \in V(G(w_0))$ with the property that for every
$v \in V(\Gamma)$, the minimal path in $G(w_0)$ from $w$ to $w_v$ does
not involve twisting at $v$.
\end{notn}

\begin{rem} $\bar{G}(w_0)$ does not in general have to contain all of the
$w_v$. 
It can also be a single vertex, as for instance in the 
two-component case if we use a multidegree which is simultaneously
concentrated on both components. However, we show below that it is always
nonempty, at least under a very mild nonnegativity hypothesis.
\end{rem}

\begin{lem}\label{lem:twist-section} Given $w \in V(G(w_0))$, and a
nonempty $S \subseteq V(\Gamma)$, suppose
that there exists a line bundle $\sL_w$ on $\widetilde{X}_0$ of multidegree 
$w$, and $s \in \Gamma(\widetilde{X}_0,\sL_w)$ such that $S$ is equal to
the set of vertices $v$ with $s|_{Z_v} \neq 0$. Then for every 
$v \in V(\Gamma)\smallsetminus S$, the set of vertices occurring as
twists in a minimal path from $w$ to $w_v$ contains $S$. For every
$v \in S$, the set of vertices occurring as
twists in a minimal path from $w$ to $w_v$ either contains $S$ or does
not contain $v$.

In particular, if $w$ is in $\bar{G}(w_0)$, then so is the twist of $w$
at all the vertices in $S$.
\end{lem}

\begin{proof} Whether or not $v \in S$, if the minimal path does not contain 
$S$, then $s$ has nonzero image under the map
$\sL_w \to \sL_{w_v}$, yielding a nonzero section of $\sL_{w_v}$. According
to Proposition 3.3 of \cite{os25}, this section cannot vanish on
$Z_v$, so we see that we must have first that $v \in S$, and second that 
$v$ does not occur in the minimal path, as desired.

The second statement follows, since if $w'$ is the twist of $w$ at the
vertices in $S$, then for any $v$ we have that the minimal path from
$w'$ to $w_v$ is obtained from the minimal path from $w$ to $w_v$ by
removing $S$, if the latter contains $S$, or by adding the complement of
$S$ otherwise.
\end{proof}

\begin{cor}\label{cor:barg} The graph $\bar{G}(w_0)$ is finite. If
there exists any $w \in V(G(w_0))$ which is everywhere nonnegative,
then $\bar{G}(w_0)$ is also nonempty, and in fact contains $w$.
\end{cor}

\begin{proof} From the definition, it is clear that if $w \in V(\bar{G}(w_0))$,
then for each $v \in V(\Gamma)$ the degree of $w$ at $v$ is bounded by the
degree of $w_v$ at $v$. Since $w$ has total degree $d$, the number of 
possible $w$ is finite. 

Now, suppose that $w \in V(G(w_0))$ is everywhere nonnegative. It is then 
clear that there exists an $\sL_w$ of multidegree $w$ (and hence an $\sL$ 
of multidegree $w_0$ having $\sL_w$ as its multidegree-$w$ twist) together 
with a section $s \in \Gamma(\widetilde{X}_0,\sL_w)$ which is nonzero on 
every component of $\widetilde{X}_0$. Then according to the first part of
Lemma \ref{lem:twist-section} with $S=V(\Gamma)$, we have that 
$w \in V(\bar{G}(w_0))$, as claimed. 
\end{proof}

\begin{rem} Without any nonnegativity condition, we may have
$\bar{G}(w_0)$ empty. However, our nonnegativity condition is quite mild:
in particular, it is implied whenever $w_0$ supports a limit $\fg^r_d$
for any $r \geq 0$. Indeed, if we have a limit $\fg^r_d$ with underlying
line bundle $\sL$, then we must have 
$\Gamma(\widetilde{X}_0,\sL_w) \neq 0$ for all $w$. In particular,
if $w$ is concentrated at $v$ for some $v$, then $w$ must have nonnegative
degree at $v$. But according to Corollary \ref{cor:concen-canon} below, 
there always exists $w$ which is concentrated at $v$ and nonnegative 
elsewhere, and we thus see that this $w$ must be nonnegative everywhere.
\end{rem}

The following corollary will not be used in the remainder of this paper,
but it is a very natural application of the preceding results.

\begin{cor} In the definition of limit linear series, it would be 
equivalent to consider \eqref{eq:gluing-map-1} only for 
$w \in V(\bar{G}(w_0))$.
\end{cor}

\begin{proof} Suppose that the desired condition is satisfied for all
$w' \in V(\bar{G}(w_0))$, so that we want to show it is also satisfied for
all other $w \in V(G(w_0))$. Given $w$, fix $w' \in V(\bar{G}(w_0))$
admitting the smallest possible path to $w$. We claim that the map
$$\Gamma(\widetilde{X}_0,\sL_{w'}) \to \Gamma(\widetilde{X}_0,\sL_{w})$$
is injective. Indeed, if $s \in \Gamma(\widetilde{X}_0,\sL_{w'})$
maps to zero, let $S \subseteq V(\Gamma)$ consist of $v$ with
$s|_{Z_v} \neq 0$, so that necessarily we have that the minimal path from
$w'$ to $w$ involves twisting at every $v \in S$. 
Now, according to Lemma \ref{lem:twist-section}, if $w''$ is obtained
from $w'$ by twisting at the vertices in $S$, then we also have
$w'' \in V(\bar{G}(w_0))$, but then the minimal path from $w''$ to $w$
is obtained by removing $S$ from minimal path from $w'$ to $w$. By
minimality, we conclude that $S=\emptyset$ and $s=0$, as claimed.
But considering each $v$ separately, we see the kernel of 
\eqref{eq:gluing-map-1} in multidegree $w'$ maps
into the kernel in multidegree $w$, so we obtain the desired statement.
\end{proof}

\section{Multidegrees and divisors on graphs}\label{sec:multidegs}

In this section, we discuss the relationship between our multidegrees
and the theory of divisors on (non-metric) graphs as developed by Baker
and Norine \cite{ba2}, \cite{b-n1}.
Some technical issues arise because
of our restriction to admissible multidegrees, but our main focus is to
make precise the close relationship between concentrated multidegrees and
$v$-reduced divisors. As an application, we prove the following theorem.

\begin{thm}\label{thm:riemann} If $r \geq g$ or $d>2g-2$, there is no
limit $\fg^r_d$ with $r>d-g$ on any $X_0$ of genus $g$.
\end{thm}

The proof given by 
Eisenbud and Harris in the compact type case is rather \textit{ad hoc},
while the proof by Amini and Baker uses their Riemann-Roch theorem for
divisors on metrized complexes. On the other hand, our proof relies on
the Riemann-Roch theorem for reducible curves. The key ingredient in
our proof is the new definition of limit linear series we have provided,
which directly relates limit linear series to dimensions of spaces of
global sections of twists of a line bundle. However, the difficulty 
remains that on reducible curves, line bundles with negative degree may
still have nonzero global sections, so we need to show that there 
always exist twists of a given line bundle for which the space of global
sections satisfies the usual bounds. We do this in Proposition 
\ref{prop:riemann} below, using the existence of $v$-reduced divisors on 
graphs. The theory of algebraic rank of divisors on graphs would give
an alternative approach; see Remark \ref{rem:alg-rank}.

First recall that a \textbf{divisor} on a graph $G$ is simply an 
integer vertex weighting (i.e., what we call a multidegree), and 
divisors $D,D'$ are \textbf{linearly equivalent} if $D'$ can be obtained 
from $D$ by a sequence of ``chip-firing moves'' at vertices of $G$, which
are exactly what we have called twists at $v$. See \S 1.3 and Lemma 4.3
of \cite{b-n1}. The following definition is the key new input from
the theory of divisors on graphs (see \S 3.1 of \cite{b-n1}):

\begin{defn}\label{defn:vred}
Given a graph $G$ and a vertex $v_0 \in V(G)$, a divisor $D$ on 
$G$ is \textbf{$v_0$-reduced} if:
\begin{nlist}
\itm $D(v) \geq 0$ for all $v \neq v_0$;
\itm for every nonempty subset $S \subseteq V(G) \smallsetminus \{v_0\}$, 
there is some $v \in S$ such that $D(v)$ is strictly smaller than the
number of edges from $v$ to $V(G) \smallsetminus S$.
\end{nlist}
\end{defn}

We recall Proposition 3.1 of \cite{b-n1}, which states:

\begin{prop}\label{prop:v-reduced} Given a graph $G$, a divisor $D$ on
$G$, and a vertex $v_0 \in V(G)$, there exists a unique divisor $D'$ on
$G$ which is linearly equivalent to $D$ and which is $v_0$-reduced.
\end{prop}

We will now develop the precise relationship between our concentrated 
(admissible) multidegrees and $v$-reducedness. A preliminary fact is the
following.

\begin{prop}\label{prop:concen-equivs} An admissible multidegree $(w,\mu)$
on $(\Gamma,\bn)$ is concentrated at $v_0 \in V(\Gamma)$ if and only if 
the induced multidegree $\widetilde{w}$ on $\widetilde{\Gamma}$ is 
concentrated at $v_0$.
\end{prop}

\begin{proof} First suppose that $(w,\mu)$ is concentrated at $v_0$,
and let $v_0,v_1,\dots$ be the ordering of $V(\Gamma)$ given by the
definition of concentrated. We extend this to an ordering on 
$V(\widetilde{\Gamma})$ as follows: between each $v_{i-1}$ and $v_i$,
for each edge $e \in E(\Gamma)$ connecting $v_i$ to some $v_{i'}$ with 
$i'<i$, if $\mu(e)=0$ we add the inserted vertices over $e$, in the order
of going from $v_{i'}$ to $v_i$. After $v_{|V(\Gamma)|-1}$, we add in
all the remaining inserted vertices, which each necessarily lie over some 
$e$ with $\mu(e)\neq 0$. In this case, for each such $e$, if the adjacent
vertices of $\Gamma$ are $v,v'$, we first add all the degree-$0$ inserted
vertices over $e$, going from $v$ to the degree-$1$ vertex, and from $v'$
to the degree-$1$ vertex, and then we finally add the degree-$1$ vertex.
It is routine to verify that this extended ordering satisfies the 
condition for $v_0$-concentratedness on $\widetilde{\Gamma}$.

Conversely, suppose that $\widetilde{w}$ is concentrated at $v_0$, and
let $v_0,v_1,\dots$ be the ordering on $V(\Gamma)$ induced by the
hypothesized ordering on $V(\widetilde{\Gamma})$; we claim that this
must have the desired property. The main observation is that in the
original ordering, if we have an edge $e$ with $\mu(e) \neq 0$, and
$v$ is adjacent to $e$, then the inserted vertex lying over $e$ and
adjacent to $v$ can't occur in the ordering until after $v$ itself has.
It then follows that negative twists
at $v_0,\dots,v_{i-1}$ on $\Gamma$ creates degree on $v_i$ at most
equal to the degree on $v_i$ obtained by negative twists on 
$\widetilde{\Gamma}$ at all the preceding vertices in the original order. 
Thus, $(w,\mu)$ is concentrated at $v_0$, as desired.
\end{proof}

\begin{prop}\label{prop:vred-conc} A multidegree on a graph $\Gamma'$ is
concentrated at $v_0$ if and only if it, considered as a divisor on
$\Gamma'$, satisfies condition (2) of Definition \ref{defn:vred}.
\end{prop}

\begin{proof} First suppose that the multidegree $w$ is concentrated at 
$v_0$, and let $v_0,v_1,\dots$ be the ordering of $V(\Gamma')$ given by the
definition. Let $S \subseteq V(\Gamma')\smallsetminus \{v_0\}$
be nonempty, and set $v=v_i$ with $i$ minimal so that $v_i$ is in $S$. 
Then by hypothesis, taking negative twists at $v_0,\dots,v_{i-1}$ results 
in the degree at $v_i$ becoming negative, which implies that the degree at
$v$ is strictly smaller than the number of edges $e$ in $E(\Gamma)$ from 
$v_i$ to $V(\Gamma) \smallsetminus S$. Thus, in this
case we find that $w$ has the desired property.

Conversely, if condition (2) of Definition \ref{defn:vred} is satisfied,
we construct the desired ordering inductively. If we have already 
found $v_0,\dots,v_{i-1}$, let 
$S := V(\Gamma')\smallsetminus \{v_0,\dots,v_{i-1}\}$, and let $v \in S$
satisfy that the degree of $w$ at $v$ is strictly less than the number
of edges from $v$ to $V(\Gamma')\smallsetminus S=\{v_0,\dots,v_{i-1}\}$.
Then we have that taking the negative twist of $w$ at $v_0,\dots,v_{i-1}$ 
makes the degree at $v$ negative, so setting $v=v_i$ produces the desired
behavior.
\end{proof}

From Propositions \ref{prop:concen-equivs} and \ref{prop:vred-conc}
we immediately conclude:

\begin{cor}\label{cor:vred-conc} An admissible multidegree $(w,\mu)$
on $(\Gamma,\bn)$ is concentrated at $v_0 \in V(\Gamma)$ if and only if 
the associated divisor $D_w$ on $\widetilde{\Gamma}$ satisfies condition
(2) of Definition \ref{defn:vred}.
\end{cor}

We also have the following, which says that, for admissible multidegrees, 
being related by twists in our sense is equivalent to the associated
divisors being linearly equivalent.

\begin{prop}\label{prop:twists-compatible} Suppose two admissible 
multidegrees $(w,\mu)$, $(w',\mu')$ have that their associated multidegrees
on $\widetilde{\Gamma}$ are related by twists. Then $(w,\mu)$, $(w',\mu')$ 
are themselves related by twists on $\Gamma$.
\end{prop}

\begin{proof} Let $\widetilde{w},\widetilde{w}'$ be the associated 
multidegrees on $\widetilde{\Gamma}$, and suppose that we can go from
$\widetilde{w}$ to $\widetilde{w}'$ by twisting $c_v$ times at $v$ for
each $v \in V(\widetilde{\Gamma})$, with each $c_v \geq 0$, and not all
$c_v>0$. Then let $(w'',\mu'')$ be the admissible 
multidegree on $\Gamma$ obtained by twisting $c_v$ times at $v$ for each 
$v \in V(\Gamma)$. If $\widetilde{w}''$ is the associated multidegree
on $\widetilde{\Gamma}$, we then have that $\widetilde{w}'$ and
$\widetilde{w}''$ are both admissible, and we can go from the latter to
the former by twisting entirely at vertices in 
$V(\widetilde{\Gamma})\smallsetminus V(\Gamma)$ (allowing negative twists). 
We claim that this implies that $\widetilde{w}'=\widetilde{w}''$. Indeed, 
this can be checked separately on each chain of inserted vertices, and
the main point is that any subchain of such chains can have its total
degree change by at most $1$ when going between admissible multidegrees
via twists on $\widetilde{\Gamma}$.
If we have any positive twists on a given chain, then the minimal
subchain containing all the vertices with positive twists will have to drop
its total degree by at least $2$ under the twists, leading to a 
contradiction. Similarly, any negative twists lead to a subchain with
total degree increasing by at least $2$. We conclude that 
$\widetilde{w}'=\widetilde{w}''$, as desired.
\end{proof}

\begin{prop}\label{prop:vred-admiss} If a multidegree $\widetilde{w}$ on 
$\widetilde{\Gamma}$ is $v$-reduced for some $v \in V(\Gamma)$, then 
it is admissible.
\end{prop}

\begin{proof} If $v_1,\dots,v_n$ is a chain of inserted vertices over an
edge $e$ of $\Gamma$, since $v \not \in \{v_1,\dots,v_n\}$, we have that
$\widetilde{w}$ is nonnegative on the $v_i$ by the definition of $v$-reduced. 
If $\widetilde{w}$ has degree at least $2$ on one of the $v_i$, setting 
$S=\{v_i\}$ violates the definition of $v$-reducedness. On the other hand, 
if $v_i,v_{i'}$ both have degree $1$ with $i<i'$, then setting 
$S=\{v_i,v_{i+1},\dots,v_{i'}\}$ likewise violates the definition.
Thus, $\widetilde{w}$ is admissible.
\end{proof}

Putting together Corollary \ref{cor:vred-conc} and Propositions
\ref{prop:twists-compatible} and \ref{prop:vred-admiss}, we thus obtain 
canonical concentrated multidegrees from the theory of $v$-reduced divisors.

\begin{cor}\label{cor:concen-canon} If $(w_0,\mu_0)$ is an admissible 
multidegree on $(\Gamma,\bn)$, for each $v \in V(\Gamma)$ there is a unique 
twist of $(w_0,\mu_0)$ which is concentrated on $v$ and nonnegative on
all $v' \neq v$.
\end{cor}

\begin{rem}\label{rem:v-reduced-role} While the canonical divisors obtained
in Corollary \ref{cor:concen-canon} are appealing and may well be important
for construction of proper moduli spaces, there are also circumstances 
where it may be better to consider other choices of concentrated 
multidegrees. For instance, if one uses the canonical multidegrees,
the finite graph $\bar{G}(w_0)$ constructed in Notation \ref{notn:barg}
below will typically have asymmetries reflecting any asymmetries in $w_0$,
while allowing degrees to go negative on some components can produce more
symmetry, and simplify resulting formulas. In addition, the canonical 
multidegrees of Corollary \ref{cor:concen-canon} need not remain 
concentrated under restriction to subcurves, which could also cause 
issues for certain arguments.
\end{rem}

In a different direction, we can also find twists of any given line bundle
on $\widetilde{X}_0$ satisfying the usual bounds for the dimension of the
space of global sections.

\begin{prop}\label{prop:riemann} Let $\sL$ be a line bundle on
$\widetilde{X}_0$ of (admissible) multidegree $w_0$, and total degree $d$.
Then there is a twist $w$ of $w_0$ such that
$$\dim \Gamma(\widetilde{X}_0,\sL_w) \leq \max(d+1-g,g).$$ 
\end{prop}

\begin{proof} Let $w_{\can}$ be the multidegree associated to the
dualizing sheaf $\omega_{\widetilde{X}_0}$ (note that this has degree $0$ 
on all inserted vertices). Choose any $v_0 \in V(\Gamma)$,
and let $w'$ be the $v_0$-reduced divisor on $\widetilde{\Gamma}$
associated to $w_{\can}-w_0$. Write $d_0$ for the degree of $w'$ on $v_0$,
so that $d_0 \leq 2g-2-d$ by the nonnegativity condition in the definition
of $v_0$-reduced. Then because $w'$ is concentrated at $v_0$, we have 
by Proposition 3.3 of \cite{os25} that the restriction map 
$$\Gamma(\widetilde{X}_0,(\omega_{\widetilde{X}_0} \otimes \sL^{-1})_{w'}) 
\to 
\Gamma(Z_{v_0},(\omega_{\widetilde{X}_0} \otimes \sL^{-1})_{w'}|_{Z_{v_0}})$$
is injective, and it follows that 
$$\dim \Gamma(\widetilde{X}_0,(\omega_{\widetilde{X}_0} \otimes \sL^{-1})_{w'}) 
\leq \max\{0,d_0+1\} \leq \max\{0,2g-1-d\}.$$
Now, $w'$ is admissible, but $w_{\can}-w'$ is not, so we modify $w'$ as
follows. For each edge $e$ of $\Gamma$, if $w'$ has degree $1$ on
an inserted vertex over $e$, then we can twist on $\widetilde{\Gamma}$
so that the degree of the two vertices of $\Gamma$ adjacent to $e$ are each 
increased by $1$, and we have degree $-1$ on some (possibly different) 
inserted vertex over $e$, and degree $0$ on the others. Apply this operation 
to each edge of $\Gamma$, and denote the resulting multidegree by $w''$.
Then we have by construction that $w_{\can}-w''$ is admissible, and is
obtained from $w_0$ by twists. Furthermore, we see that although we have
increased the degree on some vertices of $\Gamma$, this is precisely offset
by new vanishing conditions for global sections coming from chains of 
inserted vertices with negative degrees. More precisely, we still have
that $w''$ is concentrated at $v_0$, and that the image in
$\Gamma(Z_{v_0},(\omega_{\widetilde{X}_0} \otimes \sL^{-1})_{w''}|_{Z_{v_0}})$
of 
$\Gamma(\widetilde{X}_0,(\omega_{\widetilde{X}_0} \otimes \sL^{-1})_{w''})$
is canonically isomorphic to the image in
$\Gamma(Z_{v_0},(\omega_{\widetilde{X}_0} \otimes \sL^{-1})_{w'}|_{Z_{v_0}})$
of
$\Gamma(Z_{v_0},(\omega_{\widetilde{X}_0} \otimes \sL^{-1})_{w'}|_{Z_{v_0}})$,
so we still conclude that 
$$\dim 
\Gamma(\widetilde{X}_0,(\omega_{\widetilde{X}_0} \otimes \sL^{-1})_{w''}) 
\leq \max\{0,2g-1-d\}.$$
Then the Riemann-Roch theorem implies that
$$\Gamma(\widetilde{X}_0,\sL_{w_{\can}-w''})\leq d+1-g+\max\{0,2g-1-d\}=
\max\{d+1-g,g\},$$ 
as desired.
\end{proof}

We can now conclude our version of ``Riemann's theorem'' for limit linear
series.

\begin{proof}[Proof of Theorem \ref{thm:riemann}] 
The definition of a limit $\fg^r_d$ requires in particular
that we have an admissible multidegree $w_0$ and a line bundle $\sL$ of 
multidegree $w_0$ on $\widetilde{X}_0$, such 
that every twist $\sL_w$ of $\sL$ has at least an $(r+1)$-dimensional space 
of global sections. The desired statement then follows from Proposition
\ref{prop:riemann}.
\end{proof}

\begin{rem}\label{rem:alg-rank} Given our definition of limit linear series,
the theory of algebraic rank of divisors on graphs developed by Caporaso, 
Len and Melo in \cite{c-l-m1} provides a very natural alternative approach
to proving blanket non-existence results 
such as Theorem \ref{thm:riemann}, and the Clifford inequality given in
Corollary \ref{cor:riemann-clifford} below.
Indeed, to say that the algebraic
rank of a multidegree $w$ on a graph $\Gamma$ is bounded by $r$ means that 
there is some twist $w'$ of $w$ such that every line bundle of multidegree
$w'$ on every curve having dual graph $\Gamma$ has at most an 
$(r+1)$-dimensional space of global sections. In \cite{c-l-m1}, the
authors show that algebraic rank satisfies both the Riemann and
Clifford inequalities (rather like our approach, they prove the former via 
a Riemann-Roch theorem, and the latter via a comparison to ranks of 
divisors on graphs). In order to conclude corresponding bounds on limit 
linear series, one would have to check that one can replace $w'$ by an 
admissible multidegree without increasing the space of global sections.
Because the present paper already involves so many different definitions,
and we already can prove both inequalities, we do not pursue this direction.
\end{rem}

\section{Comparison to the Amini-Baker construction}\label{sec:forget}

In this section, we show that in full generality, there is a forgetful
map from our limit linear series to Amini-Baker limit linear series.
Accordingly, we begin by recalling the definitions of Amini and Baker.
Because we have (following \cite{os25}) already used $\Gamma$ and $G$
in our setup, we will instead use $\GamAB$ and $\GAB$ for the Amini-Baker 
usage of a metric graph with an imbedded finite graph. Before recalling
the Amini-Baker definition, we recall the corresponding definitions for
metric graphs, since we will need them as well.

\begin{defn}\label{def:trop-BN} Let $\GamAB$ be a metric graph.
A \textbf{divisor} $D$ on $\GamAB$ is a 
finite formal $\ZZ$-linear sum $\sum_i a_i [x_i]$, where each $x_i$ is a
point of $\GamAB$. The \textbf{degree} of $D$
is $\sum_i a_i$, and $D$ is \textbf{effective} if $a_i \geq 0$ for all $i$.

A \textbf{nonzero rational function} $f$ on $\GamAB$ is a 
(continuous) piecewise linear function on $\GamAB$, with each piece
having integer slope. The divisor $\dv f$ associated to $f$ is
defined in terms of slopes as follows:
for each $x \in \GamAB$, the coefficient of $[x]$ in
$\dv f$ is the sum of the outgoing slopes of $f$ at $x$.

Two divisors on $\GamAB$ are \textbf{linearly equivalent} if their difference
is $\dv f$ for some nonzero rational function $f$ on $\GamAB$.

A divisor $D$ on $\GamAB$ has \textbf{rank} $r$ if $r$ is maximal such that 
for all effective divisors $E$ on $\GamAB$ of degree $r$, we have that
$D-E$ is linearly equivalent to an effective divisor.
\end{defn}

We now recall the Amini-Baker definitions.

\begin{defn} A \textbf{metrized complex of curves} $\fC$ consists of a 
connected finite loopless 
graph $\GAB$ together with a length function on $E(\GAB)$, 
a smooth projective curve $C_v$ associated to each $v \in V(\GAB)$, and
for each $v \in V(\GAB)$, a bijection between the set of $e \in E(\GAB)$ which
are adjacent to $v$, and a subset $\cA_v=\{x^e_v\}$ of points of $C_v$. 
We will denote by $\GamAB$
the metric graph induced by $\GAB$ together with the edge weights.
\end{defn}

We will use $C_v$ and $Z_v$ relatively interchangeably in the following,
although when we are unambiguously in the Amini-Baker context or in our
own context we will generally use $C_v$ and $Z_v$ respectively, and we will
use $Z_v$ when we want to think of it as a component of $X_0$ (for instance,
if we need to refer to the smooth locus of $Z_v$).

\begin{rem}\label{rem:a-b-loops} Amini and Baker allow their graphs to
have loops. However, in \cite{os25} for the sake of simplicity components
were not allowed to have self-nodes, so we will assume throughout that
our graphs do not have loops.
\end{rem}

\begin{defn} A \textbf{divisor} $\cD$ on a metrized complex of curves
$\fC$ is a finite formal $\ZZ$-linear sum
$\sum_i a_i [x_i]$, where each $x_i$ is either a point of some $C_v$,
or a point of $\GamAB\smallsetminus V(\GAB)$. The \textbf{degree} of $\cD$
is $\sum_i a_i$, and $\cD$ is \textbf{effective} if all $a_i$ are
nonnegative. We denote by $D_v$ the part of $\cD$ supported on
$C_v$, and by $D_{\GamAB}$ the divisor on $\GamAB$ obtained as
$$\sum_{v \in V(\GAB)} \deg D_v [v]
+\sum_{i:x_i \in \GamAB\smallsetminus V(\GAB)} a_i x_i.$$
Thus, $\deg D_{\GamAB}=\deg \cD$.

A \textbf{nonzero rational function} $\ff$ on $\fC$ consists of a 
nonzero rational function $f_v$ on each $C_v$, and a nonzero
rational function $f_{\GamAB}$ on $\GamAB$.
The divisor $\dv \ff$ associated to $\ff$ is
$\dv f_{\GamAB}+\sum_v \dv f_v$, where $\dv f_v$ is the usual divisor
on $C_v$, and $\dv f_{\GamAB}$ is defined in terms of slopes as follows:
for $x \in \GamAB\smallsetminus V(\GAB)$, the coefficient of $[x]$ in
$\dv f_{\GamAB}$ is the sum of the outgoing slopes of $f_{\GamAB}$ at $x$,
while for $x^e_v \in \fC$, the coefficient of $[x^e_v]$ in 
$\dv f_{\GamAB}$ is the outgoing slope of $f_{\GamAB}$ at $v$ in the
direction of $e$.

Two divisors on $\fC$ are \textbf{linearly equivalent} if their difference
is $\dv \ff$ for some nonzero rational function $\ff$ on $\fC$.

A divisor $\cD$ on $\fC$ has \textbf{rank} $r$ if $r$ is maximal such that 
for all effective divisors $\cE$ on $\fC$ of degree $r$, we have that
$\cD-\cE$ is linearly equivalent to an effective divisor. If we are also
given a tuple $(H_v)_v$ with each $H_v$ an $(r+1)$-dimensional
subspace of the function field $K(C_v)$, the \textbf{restricted rank}
of $\cD$ with respect to the $(H_v)_v$ is the maximal $r'$ such that
for all effective divisors $\cE$ on $\fC$ of degree $r'$, we can find
a nonzero rational function $\ff$ on $\fC$ such that each $f_v$ is in
$H_v$, and such that $\cD-\cE-\dv \ff$ is effective.
\end{defn}

Note that in the above situation, the restricted rank of $\cD$ with
respect to $(H_v)_v$ is always at most $r$.
Also note that in the context of metrized complexes, $\dv f_{\GamAB}$ is not a 
divisor on $\GamAB$, because what would
be its support at $v \in V(G)$ is divided among the $x^e_v$. However, we do
have the compatibility that the divisor $(\dv f_{\GamAB})_{\GamAB}$ on 
$\GamAB$ is equal to the usual divisor on $\GamAB$ constructed from
$f_{\GamAB}$ in the context of metric graphs. We will attempt to state
clearly which form of $\dv f_{\GamAB}$ is being used when.

The Amini-Baker definition of limit linear series is then as follows.

\begin{defn}\label{def:ab-lls}
A \textbf{limit} $\fg^r_d$ on a metrized complex of curves
$\fC$ is an equivalence class of pairs $(\cD,(H_v)_v)$, where $\cD$ is
a divisor of degree $d$ on $\fC$, each $H_v$ is an $(r+1)$-dimensional
subspace of $K(C_v)$, and the restricted rank of $\cD$ with respect to
$(H_v)_v$ is equal to $r$. The equivalence relation is given by
nonzero rational functions $\ff$ by replacing $\cD$ with $\cD+\dv \ff$ and 
each $H_v$ with $H_v/f_v$.
\end{defn}

We now connect our setup with that of Amini and Baker as follows.

\begin{notn}\label{notn:assoc-complex} Given a nodal curve $X_0$ with
dual graph $\Gamma$, and a chain structure $\bn$ on $\Gamma$, let 
$\fC(X_0,\bn)$ denote the metrized complex of curves obtained by setting
$\GAB=\Gamma$, $\GamAB$ equal to the metric graph obtained from $\Gamma$
by letting $\bn$ specify the lengths of each edge, and for $v \in V(\Gamma)$,
setting $C_v$ to be the corresponding component of $X_0$.
\end{notn}

To describe the construction relating our limit linear series to
Amini-Baker limit linear series, we first need the following definition,
which is entirely in the former context. It generalizes the divisor 
sequences constructed
in the pseudocompact type case in Notation 5.8 of \cite{os25}.

\begin{notn}\label{notn:Dwv}
Given $w \in V(G(w_0))$ and $v \in V(\Gamma)$, let 
$D_{w,v}$ be the divisor on $Z_v$ defined as follows:
if 
$P(w,v_1,\dots,v_m)=(w_1,\mu_1),(w_2,\mu_2),\dots,(w_{m+1},\mu_{m+1})$ 
is a path 
from $w$ to $w_v$, so that $w=(w_1,\mu_1)$ and $w_v=(w_{m+1},\mu_{m+1})$, 
let $S \subseteq \{1,\dots,m\}$ consist of $i$ such
that $v_i$ is adjacent to $v$, and for $i \in S$, let $e_i\in E(\bar{\Gamma})$ 
be the connecting edge. Then set
$$D_{w,v}=\sum_{i \in S} \ \sum_{\substack{\tilde{e}\text{ over }e_i:\\
\mu_{i+1}(\tilde{e})=0}} \!\!\! P_{\tilde{e}}\ -
\sum_{i:v_i =v} \ \sum_{\substack{\tilde{e}\text{ adjacent to }v:\\
\mu_{i}(\tilde{e})=0}} \!\!\!\!\!\! P_{\tilde{e}}.$$
\end{notn}

The point of this definition is the following.

\begin{prop}\label{prop:twist-divs} In the situation of Notation
\ref{notn:Dwv}, we have that 
$$\sL_w|_{Z_v} \cong \sL^v(-D_{w,v}).$$ 
If further $w \in V(\bar{G}(w_0))$, then $D_{w,v}$ is effective, and
the restriction to $Z_v$ of any section of 
$\sL_w$ will be contained in $\Gamma(Z_v,\sL^v(-D_{w,v}))$.
\end{prop}

The verification of this is routine. We also have the following.

\begin{prop}\label{prop:dwv-indep} $D_{w,v}$ is independent of the choice
of path from $w$ to $w_v$.
\end{prop}

\begin{proof} We first verify that $D_{w,v}$ is independent of reordering.
It is enough to consider
iterated swapping of $v_i$ with $v_{i+1}$. Because $\mu_i(\tilde{e})$ is
only changed by twisting at $v'$ if $v'$ is adjacent to $\tilde{e}$, it is 
clear from the definition 
that such a swap can only affect $D_{w,v}$ when either $i \in S$ and
$v_{i+1}=v$ or vice versa. Moreover, only the $P_{\tilde{e}}$ terms with
$\tilde{e}$ adjacent to $v_i$ and $v_{i+1}$ can be affected by such a swap. 
In either of these cases, the individual terms appearing
in the formula for $D_{w,v}$ may change for such $P_{\tilde{e}}$,
but we see that for such $\tilde{e}$, the $P_{\tilde{e}}$ terms contributed
by $v_i$ and $v_{i+1}$ always cancel if they are nonzero, so the swap
can't affect the overall formula for $D_{w,v}$.

To conclude independence of path, it suffices to verify that $D_{w,v}$ is
also not affected by appending one of each vertex of $\Gamma$ to any
given path. Given independence of ordering, we may assume that $v$ is
appended last, and in this case it is clear that the negative contribution
from the twist at $v$ precisely cancels the positive contributions from the
prior twists at the other vertices, so $D_{w,v}$ is again unchanged.
\end{proof}

\begin{notn} Given $(w,\mu) \in V(G(w_0))$, let $D_{\mu}$ be the divisor
on $\GamAB$ obtained by summing over, for each $e \in E(\Gamma)$ with 
$\mu(e) \neq 0$, the point on $e$ at distance $\mu(e)$ from the tail of
$e$.
\end{notn}

We now relate our limit linear series to Amini-Baker limit linear series
via the following construction:

\begin{defn}\label{def:forget} Fix $(\sL,(V^v)_{v \in V(\Gamma)})$ with
$\sL$ a line bundle of $\widetilde{X}_0$ of multidegree $w_0$, and each
$V^v$ an $(r+1)$-dimensional subspace of the resulting
$\Gamma(Z_v,\sL^v)$. 
For each $v \in V(\Gamma)$, choose a nonzero $s_v \in V^v$, and fix
also a $w \in V(G(w_0))$.
Then let $\cD$ be the divisor on $\fC(X_0,\bn)$ given by
$$\left(\sum_{v \in V(\Gamma)} \dv s_v-D_{w,v}\right) +D_{\mu}.$$
For each $v \in V(\Gamma)$, let 
$$H_v=\left\{\frac{s}{s_v}: s \in V^v\right\}.$$
\end{defn}

\begin{prop}\label{prop:forget-independence} In the construction of
Definition \ref{def:forget}, different choices of $w$, of the $s_v$,
or of the concentrated multidegrees $w_v$ will yield equivalent tuples
$(\cD,(H_v)_v)$ in the sense of Definition \ref{def:ab-lls}.
\end{prop}

\begin{proof} First, it is clear that modifying one of the $s_v$ changes
$\cD$ and $H_v$ by the divisor of a rational function on $Z_v$.
Next, we see that if we replace $w$ by $w'$, the resulting 
$\cD$ is modified by the divisor of a rational function on $\GamAB$:
specifically, by the divisor of the $f_{\GamAB}$ such that 
$w'=w+\dv f_{\GamAB}$, where
$\dv f_{\GamAB}$ is considered as a divisor on $\widetilde{\Gamma}$.
To prove this, it is clearly enough to consider 
the case that $w'$ is obtained from $w$ by twisting at a single 
$v \in V(\Gamma)$.
In this case, for any $v' \in V(\Gamma)$ a minimal
path $P'_{v'}$ from $w'$ to $w_{v'}$ can be obtained from a (possibly 
nonmininal) path $P_{v'}$ from $w$ to $w_{v'}$ by removing a twist at $v$.
If $v'$ is neither equal to $v$ nor adjacent to $v$, it is then clear that
we have $D_{w,v'}=D_{w',v'}$. We will have $D_{w',v}$ obtained from
$D_{w,v}$ by adding $P_{\tilde{e}}$ for all $\tilde{e}$ adjacent to $v$
with $\mu(\tilde{e})=0$.
Finally, if $v'$ is adjacent to $v$, we will have
$D_{w',v'}$ obtained from
$D_{w,v'}$ by subtracting $P_{\tilde{e}}$ for all $\tilde{e}$ adjacent to $v$
with $\mu'(\tilde{e})=0$. Here $\mu$ and $\mu'$ come from $w$ and $w'$
respectively. 
It is then routine to check that the overall change
in $\cD$ is precisely what is obtained as $\dv f_{\GamAB}$ where 
$f_{\GamAB}$ comes from twisting at $v$.

Finally, suppose that we replace a concentrated multidegree $w_v$ by
some $w'_v$. We will show that in this case, in fact both $\cD$ and the
$H_v$ are unchanged. In order to show this, we first 
need to
recall the construction given in Proposition 3.5 of \cite{os25}, which
begins with the observation that it is enough to treat the case that
$w'_v$ is obtained from $w_v$ by twisting at vertices other than $v$.
In this case, we construct a $V'^v$ in multidegree $w'_v$ simply by
taking the image of $V^v$ under the natural map 
$\sL^v \to \sL'^v:=\sL_{w'_v}|_{Z_v}$, which is injective because no
twists at $v$ were required to go from $w_v$ to $w'_v$. We then see that
if we replace a given choice of $s_v$ by its image in $V'^v$, the space
$H_v$ is in fact completely unchanged under this procedure. On the
other hand, both $\dv s_v$ and $D_{w,v}$ are modified by the same effective
divisor determined by the twists to go from $w_v$ to $w'_v$, so we see
that $\cD$ is also unchanged.
We thus conclude the asserted independence of choices.
\end{proof}

Our main comparison result is the following:

\begin{thm}\label{thm:forget} The tuple $(\cD,(H_v)_{v \in V(\Gamma)})$
of Definition \ref{def:forget}
is a limit $\fg^r_d$ in the sense of Amini and Baker.
\end{thm}

\begin{proof} First, by construction $\deg(\dv s_v - D_{w,v})$ is equal
to the degree of $w$ at $v$ for each $v \in V(\Gamma)$, so $\cD$ has degree 
$d$.
We next show that the restricted rank of $\cD$ is at least (hence exactly)
$r$. According to Theorem \ref{thm:a-b-r-rk}, it is enough to consider 
effective divisors $E$ of degree $r$ supported on the smooth loci of the
$Z_v$. 
Fix such an $E$; we need to show that
$\cD-E$ is linearly equivalent to an effective divisor via the given
spaces $H_v$. According to Proposition \ref{prop:forget-independence}
(and noting that in our construction, we always have $1 \in H_v$),
it is enough to show that there exist choices of $w$ and $s_v$ such that
the resulting $\cD$ has $\cD-E$ effective. 
 For each $v \in V(\Gamma)$, denote by $E_v$ the part of
$E$ lying on $C_v$. For each $w \in V(\bar{G}(w_0))$, define 
$S_w \subseteq V(\Gamma)$ to be the set of $v$ such that 
$V^v(-D_{w,v}-E_v) \neq 0$. Thus, if we show that for some $w$ we have
$S_w= V(\Gamma)$, we conclude the desired statement on the restricted rank
of $\cD$. Now, define also
$S'_w \subseteq V(\Gamma)$ to be the set of $v$ such that there is a 
global section of $\sL_w(-E)$ which is in the kernel of
\eqref{eq:gluing-map-1} and which does not vanish identically on $Z_v$.
Note that by taking linear combinations, we can then find such a global 
section of $\sL_w(-E)$ which does not vanish identically on any $Z_v$
with $v \in S'_w$. Because $w \in \bar{G}(w_0)$, the subspace of 
$\Gamma(\widetilde{X}_0, \sL_w(-E))$ lying in the kernel of 
\eqref{eq:gluing-map-1} must have the property that any section with
nonvanishing restriction to a given $Z_v$ must have nonzero image in
$V^v(-D_{w,v}-E_v)$,
so we have $S'_w \subseteq S_w$.
Moreover, by definition of a limit linear series, the kernel of
\eqref{eq:gluing-map-1} has dimension at least $r+1$, so since $\deg E=r$, 
we have that $S'_w \neq \emptyset$. 

Starting from
any $w \in V(\bar{G}(w_0))$, if $S'_w \neq V(\Gamma)$, let $w'$ be the 
multidegree obtained from 
$w$ by twisting at all the vertices in $S'_w$. It follows from 
Lemma \ref{lem:twist-section} that $w'$ remains in $V(\bar{G}(w_0))$.
By construction, we see that we also have $S_w \subseteq S_{w'}$: indeed,
$S'_w \subseteq S_{w'}$ because the global section supported on $S'_w$
still yields nonzero elements of each $V^v(-D_{w,v}-E_v)$ for $v \in S'_w$
after twisting by the vertices in $S'_w$. On the other hand, for 
$v \not\in S'_w$, twisting by the vertices in $S'_w$ can only increase
$V^v(-D_{w,v}-E_v)$, so also $S_w \smallsetminus S'_w \subseteq S_{w'}$.
If we repeat this process, by the finiteness of $V(\bar{G}(w_0))$ we
must either eventually have $S'_w=V(\Gamma)$ and hence that $S_w=V(\Gamma)$, 
or that we return to
a $w$ which we had previously been at. In the latter case, we will
have necessarily twisted at every $v \in V(\Gamma)$, so every $v$ will
have occurred in $S'_{w'}$ for some $w'$, and it follows again that 
$S_w =V(\Gamma)$. The theorem follows.
\end{proof}

\begin{rem} While it is natural to expect that the $\cD$ of Definition
\ref{def:forget} is always obtained as the divisor of a global section
in a suitable multidegree, this is not necessarily the case. Indeed,
there exist (crude) Eisenbud-Harris limit $\fg^0_d$s on $2$-component
curves of compact type for which there is no multidegree supporting a
section which is nonzero on both components.
\end{rem}

As an application of Theorem \ref{thm:forget}, we conclude a new proof of
Theorem \ref{thm:riemann} for our limit linear series, as well as
a version of Clifford's theorem.

\begin{cor}\label{cor:riemann-clifford} If $r \geq g$ or $d>2g-2$, there 
is no limit $\fg^r_d$ with $r>d-g$ on any $X_0$ of genus $g$.

If $r \leq g$ or $d \leq 2g$, there 
is no limit $\fg^r_d$ with $2r>d$ on any $X_0$ of genus $g$.
\end{cor}

\begin{proof} Amini and Baker \cite{a-b1} prove the corresponding statements 
for their divisors as a formal consequence of their Riemann-Roch theorem 
(see their Theorem 3.2, Theorem 3.4, and Remark 5.8). More specifically,
for the latter statement, they prove that if a divisor $\cD$ is special,
then $2r \leq d$. But their Riemann-Roch theorem immediately implies that 
if $r \leq g$, then $d \leq 2g$, and also that if $\cD$ is nonspecial and 
$d \leq 2g$, the inequality $2r \leq d$ holds. Thus, the desired 
statements hold for Amini-Baker divisors, and hence for Amini-Baker 
limit linear series, and the corollary then follows from Theorem 
\ref{thm:forget}.
\end{proof}

Having produced a forgetful map on limit linear series on a nodal curve,
we should also verify that our construction behaves well in a smoothing
family. Specifically, we show that it is compatible with the previously
constructed specialization maps to limit linear series \cite{os25} and to 
Amini-Baker limit linear series \cite{a-b1}, as follows:

\begin{prop}\label{prop:specialize} Let $B$ be the spectrum of a discrete
valuation ring, and $\pi:X \to B$ flat and proper with nodal special fiber
$X_0$ having smooth components, and smooth generic fiber $X_{\eta}$. 
Let $\fC$ be the metrized complex of curves induced by $\pi$.

Then the specialization map on $\fg^r_d$s on $X_{\eta}$ constructed by 
Amini and Baker (Theorem 5.9 of \cite{a-b1}) is equal to the composition 
of the specialization map constructed in Corollary 3.15 of \cite{os25} with
the forgetful map of Definition \ref{def:forget}.
\end{prop}

\begin{proof}
Let $\bn$ be the chain structure induced on $X_0$ by the singularities of
$X$, and let $\widetilde{\pi}:\widetilde{X} \to B$ be the regularization of 
$X$, so that $\fC=\fC(X_0,\bn)$. Let $(\sL_{\eta},V_{\eta})$ be a $\fg^r_d$
on $X_{\eta}$. Let $\widetilde{\fC}$ be the metrized complex of 
$\widetilde{\pi}$, obtained from $\fC$ by placing a copy of $\PP^1_k$ at 
every integral internal point of
every edge of $\GamAB$. Let $(\cD,(H_v)_v)$ on $\fC$ and
$(\widetilde{\cD},(\widetilde{H}_{\tilde{v}})_{\tilde{v}})$ on 
$\widetilde{\fC}$ be obtained by specializing $(\sL_{\eta},V_{\eta})$.
Then Proposition \ref{prop:specialize-ab} says that 
$H_v=\widetilde{H}_v$ for all $v \in V(\Gamma)$, and that $\cD$ is
obtained from $\widetilde{\cD}$ by replacing, for each
$\tilde{v}\in V(\widetilde{\Gamma})\smallsetminus V(\Gamma)$, 
any points supported on $C_{\tilde{v}}$ with 
with the same number of points supported at the point of $\GamAB$ 
corresponding to $\tilde{v}$. On the other hand, the discussion in
\cite{a-b1} preceding Theorem 5.10 describes 
$(\widetilde{\cD},(\widetilde{H}_{\tilde{v}})_{\tilde{v}})$ as follows:
$\widetilde{\cD}$ is obtained by choosing any extension of $\sL_{\eta}$ to 
$\widetilde{X}$ and then restricting to the components of 
$\widetilde{X}_0$ and choosing any sections. For each $\tilde{v}$, we then
let $\widetilde{\cD}_{\tilde{v}}$ be the unique $\tilde{v}$-reduced divisor 
linearly equivalent to $\widetilde{\cD}$ such that 
$\widetilde{\cD}-\widetilde{\cD}_{\tilde{v}}$ is equal to the divisor of a 
rational function on $\GamAB$ (see \S 3.1 of \cite{a-b1} for the 
definition of $\tilde{v}$-reduced in the metrized complex setting).
Note that this gives an admissible multidegree on $\widetilde{X}_0$, so
let $\sL_{\tilde{v}}$ be the extension of $\sL_{\eta}$ in the multidegree 
determined by $\widetilde{\cD}_{\tilde{v}}$. Then the restriction 
$V^{\tilde{v}}$ of
$V_{\eta} \cap \Gamma(\widetilde{X},\sL_{\tilde{v}})$ to $Z_{\tilde{v}}$ is
an $(r+1)$-dimensional space of sections,
and if $D^{\tilde{v}}_{\tilde{v}}$ denotes the $C_{\tilde{v}}$-part of 
$\widetilde{\cD}_{\tilde{v}}$, the sections in $V^{\tilde{v}}$ have
divisors linearly equivalent
to $D^{\tilde{v}}_{\tilde{v}}$, so we can set
$\widetilde{H}_{\tilde{v}}\subseteq 
\Gamma(C_{\tilde{v}},\sO_{C_{\tilde{v}}}(D^{\tilde{v}}_{\tilde{v}}))
\subseteq K(C_v)$
be the resulting $(r+1)$-dimensional subspace induced by $V^{\tilde{v}}$.

Then by construction we also have that the divisors of sections of 
$V^{\tilde{v}}$
all occur as $D^{\tilde{v}}_{\tilde{v}}+f_v$ for some $f_v \in
\widetilde{H}_{\tilde{v}}$. Thus, if the multidegrees induced by the
$\cD^{\tilde{v}}$ are concentrated at $\tilde{v}$, this construction 
agrees with our specialization construction combined with Definition
\ref{def:forget}. 
On the other hand, even if the multidegrees are not
concentrated, by definition of $\tilde{v}$-reducedness on a metrized 
complex, their behavior under restriction to $C_{\tilde{v}}$ for the
line bundles in question are the same as the concentrated case. Thus, we
can apply the same argument as in the proof of independence of the
choice of the $w_v$ in Proposition \ref{prop:forget-independence}
to conclude the desired compatibility.
\end{proof}

\begin{rem}\label{rem:clifford} One might wonder whether one can
prove Clifford's inequality via the Riemann-Roch theorem for reducible
curves, as in the proof of Theorem \ref{thm:riemann}, without making use
of the Amini-Baker theory. This is presumably possible, but the main 
difficulty which has to be addressed is that in order to obtain an inequality 
$$h^0(X_0,\omega_{X_0}) 
\geq h^0(X_0,\sL_w) + h^0(X_0,\omega_{X_0} \otimes \sL_w^{-1}),$$
one needs to control which components of $X_0$ the relevant global sections
vanish on. This suggests that such an argument would necessarily involve 
combinatorial ideas rather similar to ranks of divisors on graphs and/or 
complexes of curves. From this point of view, the proof we have given 
appears quite natural.
\end{rem}

\begin{rem}\label{rem:v-reduced-vers} Note that the Amini-Baker definition
of $v$-reducedness takes into account the geometry of each $D_{v'}$ on 
$C_{v'}$. Thus 
it is not purely combinatorial, and in particular, the induced multidegree
need not be concentrated in our sense (equivalently, the induced divisor
on the underlying graph need not be $v$-reduced). However, for the 
particular $\cD$ in question, if it is $v$-reduced in the Amini-Baker
sense it still behaves under restriction to $C_v$ as if it were
concentrated at $v$.
\end{rem}

\section{The pseudocompact type case}\label{sec:forget-pct}

We now suppose that our curve $X_0$ is of pseudocompact type. In this 
case, we have an equivalent characterization of our limit linear
series, which depends on some additional notation.

\begin{defn}\label{def:pseudocompact}
If $\Gamma$ is a graph (possibly with multiple edges), let $\bar{\Gamma}$
be the graph (without multiple edges) having the same vertex set as 
$\Gamma$, and with a single edge between any pair of vertices which are 
adjacent in $\Gamma$. We say that $\Gamma$ is a \textbf{multitree} if
$\bar{\Gamma}$ is a tree. We say that a nodal curve is of 
\textbf{pseudocompact type} if its dual graph is a multitree.
\end{defn}

When $\Gamma$ is a
multitree, we also have a well-defined notion of twisting on one side
of a node, as follows:

\begin{defn}\label{def:twist-node} If $\Gamma$ is a multitree,
and $(e,v)$ a pair of an edge $e$ and an adjacent 
vertex $v$ of $\bar{\Gamma}$, 
given an admissible multidegree $w$, we define the \textbf{twist} of
$w$ at $(e,v)$ to be obtained from $w$ by twisting at all $v'$ which
lie on the same connected component as $v$ in 
$\bar{\Gamma}\smallsetminus\{e\}$.
\end{defn}

This twist can be described explicitly as follows: if $v'$ is the other
vertex adjacent to $e$, everything is unchanged except at $v,v'$ and the
edges of $\Gamma$ over $e$. For each $\tilde{e}$ of $\Gamma$ over $e$, 
the twist increases $\mu(\tilde{e})$ by $\sigma(\tilde{e},v)$. It decreases 
$w_{\Gamma}(v)$ by the number of $\tilde{e}$ for which $\mu(\tilde{e})$ had 
been equal to $0$, and it increases $w_{\Gamma}(v')$ by the number of
$\tilde{e}$ for which the new $\mu(\tilde{e})$ is zero.

To simplify the situation, we assume the following:

\begin{sit}\label{sit:pct} Suppose that the dual graph $\Gamma$ of 
$X_0$ is a multitree, fix any admissible multidegree $w_0$, and let
$(w_v)_{v \in V(\Gamma)}$ be a
collection of elements of $V(G(w_0))$ such that:
\begin{Ilist}
\itm each $w_v$ is concentrated at $v$, and nonnegative away from $v$;
\itm for each $v,v' \in V(\bar{\Gamma})$ connected by an edge $e$, the
multidegree $w_{v'}$ is obtained from $w_v$ by twisting $b_{v,v'}$ times at
$(e,v)$, for some $b_{v,v'}\in \ZZ_{\geq 0}$.
\end{Ilist}
\end{sit}

Thus, according to Corollary \ref{cor:vred-conc}, the $w_v$ are simply 
the $v$-reduced
divisors on $\widetilde{\Gamma}$ linearly equivalent to $w_0$; however,
the existence of $w_v$ satisfying the above conditions is easy to see 
directly in the pseudocompact type case -- see
Proposition 2.9 of \cite{os23}.

We also have the following, which states that the conditions of Situation
5.4 of \cite{os25} are satisfied.

\begin{prop}\label{prop:pct-restrict} In Situation \ref{sit:pct},
each $w_v$ remains concentrated under restriction to a connected subcurve
containing $Z_v$.
\end{prop}

\begin{proof} Since $w_v$ is nonnegative away from $v$, the ordering
of vertices in the definition of concentratedness must be compatible
with distance from $v$ in $\bar{\Gamma}$, in the sense that a given $v'$
cannot appear until the adjacent vertex in the direction of $v$ has already 
appeared.
This implies that the concentration condition is preserved under 
restriction, using the same ordering.
\end{proof}

We now have to review some notation in order to state our equivalent 
characterization of limit linear series in the pseudocompact type case.

\begin{notn}\label{not:twist-divs} In Situation \ref{sit:pct},
for each pair $(e,v)$ of an
edge and adjacent vertex of $\bar{\Gamma}$, let 
$D^{(e,v)}_0,\dots,D^{(e,v)}_{b_{v,v'}+1}$ be the sequence of 
divisors on $Z_v$ defined by 
$$D^{(e,v)}_i:=D_{w_i,v},$$
where $w_i$ is obtained from $w_v$ by twisting $i$ times at $(e,v)$.
\end{notn}

Note that going from $w_i$ to $w_v$ does not involve twists at $v$, so
the $D^{(e,v)}_i$ are effective. In addition, the above definition agrees 
with that of Notation 5.8 of \cite{os25}. We will use the above divisor
sequences to construct generalized vanishing sequences as follows.

\begin{defn}\label{def:multivanishing} Let $X$ be a smooth projective
curve, $r,d \geq 0$, and $D_0 \leq D_1 \leq \dots \leq D_{b+1}$ a
sequence of effective divisors on $X$, with $D_0=0$ and $\deg D_{b+1}>d$.
Given $(\sL,V)$ a $\fg^r_d$ on $X$, define the
\textbf{multivanishing sequence} of $(\sL,V)$ along $D_{\bullet}$ to
be the sequence
$$a_0 \leq \dots \leq a_r$$
where a value $a$ appears in the sequence $m$ times if for some $i$ we
have $\deg D_i=a$, $\deg D_{i+1}>a$, and
$\dim \left(V(-D_i)/V(-D_{i+1})\right)=m$.

Also, given $s \in V$ nonzero, define the \textbf{order of vanishing}
$\ord_{D_{\bullet}} s$ along $D_{\bullet}$ to be $\deg D_i$, where
$i$ is maximal so that $s \in V(-D_i)$.

Finally, we say that $i$ is \textbf{critical} for $D_{\bullet}$ if
$D_{i+1} \neq D_i$.
\end{defn}

We recall (Proposition 4.6 of \cite{os25}) that we have gluing maps as
follows:

\begin{prop}\label{prop:twists-basic} In the situation of Notation 
\ref{not:twist-divs}, suppose we also have a line bundle $\sL$ on 
$\widetilde{X}_0$ with induced twists $\sL_w$ as in Notation
\ref{not:more-twist}. 
Then considering chains of exceptional curves between $Z_v$ and $Z_{v'}$ on
which $\sL_{w_i}$ is trivial, we obtain isomorphisms
$$\vp_i:\sL^{(e,v)}(-D^{(e,v)}_i)/\sL^{(e,v)}(-D^{(e,v)}_{i+1}) \risom
\sL^{(e,v')}(-D^{(e,v')}_{b_{v,v'}-i})/
\sL^{(e,v')}(-D^{(e,v')}_{b_{v,v'}+1-i})$$
for each $i$.
\end{prop}

The following is Theorem 5.9 of \cite{os25}.

\begin{thm}\label{thm:equiv} In the situation of Definition
\ref{def:lls}, suppose further that $X_0$ is of pseudocompact type,
and we are in Situation \ref{sit:pct}.
Then given a tuple $(\sL,(V^v)_{v \in V(\Gamma)})$,
for each pair $(e,v)$ in $\bar{\Gamma}$, let $a^{(e,v)}_0,\dots,a^{(e,v)}_r$
be the multivanishing sequence of $V^v$ along $D^{(e,v)}_{\bullet}$.
Then $(\sL,(V^v)_v)$ is a limit linear series if and only if
for any $e \in E(\Gamma)$, with adjacent vertices $v,v'$, we have:
\begin{Ilist}
\itm for $\ell=0,\dots,r$, if $a^{(e,v)}_{\ell}=\deg D^{(e,v)}_j$ with $j$
critical for $D^{(e,v)}_{\bullet}$, then
\begin{equation}\label{eq:eh-genl}
a^{(e,v')}_{r-\ell} \geq \deg D^{(e,v')}_{b_{v,v'}-j};
\end{equation}
\itm there exist bases $s^{(e,v)}_0,\dots,s^{(e,v)}_r$ of $V^v$ and
$s^{(e,v')}_0,\dots,s^{(e,v')}_r$ of $V^{v'}$ such that
$$\ord_{D^{(e,v)}_{\bullet}}s^{(e,v)}_{\ell} = a^{(e,v)}_{\ell},
 \quad \text{ for } \ell=0,\dots,r,$$
and similarly for $s^{(e,v')}_{\ell}$,
and for all $\ell$ with \eqref{eq:eh-genl} an equality, we have
$$\vp^{(e,v)}_{j}(s^{(e,v)}_{\ell})=s^{(e,v')}_{r-\ell}$$
when we consider $s^{(e,v)}_{\ell} \in V^v(-D^{(e,v)}_{j})$ and
$s^{(e,v')}_{r-\ell} \in V^{v'}(-D^{(e,v')}_{b_{v,v'}-j})$, where
$j$ is as in (I), and $\vp_j$ is
as in Proposition \ref{prop:twists-basic}.
\end{Ilist}
\end{thm}

We refer to condition (I) above as the \textbf{multivanishing inequality}
and condition (II) as the \textbf{gluing condition}. Our main purpose
here is to show that we still obtain a forgetful map to Amini-Baker
limit linear series if we drop the gluing condition. Under these 
circumstances, we lose control over spaces of global sections, so our
argument for Theorem \ref{thm:forget} no longer applies. Instead, we
use a different argument to conclude the following.

\begin{prop}\label{prop:forget-pseudocompact}
In the situation of Theorem \ref{thm:equiv}, suppose that 
$(\sL,(V^v)_{v \in V(\Gamma)})$ satisfies condition (I) of the theorem,
and construct $(\cD,(H_v)_{v \in V(\Gamma)})$ as in Definition 
\ref{def:forget}. Then $(\cD,(H_v)_{v \in V(\Gamma)})$
is a limit $\fg^r_d$ in the sense of Amini and Baker.
\end{prop}

\begin{proof} As in the proof of Theorem \ref{thm:forget}, we have to
show that the restricted rank of $\cD$ is at least $r$, it is enough 
to consider effective divisors $E$ of degree $r$ supported on the $C_v$,
and we want to show that there exists 
$w \in V(\bar{G}(w_0))$ such that $V^v(-D_{w,v}-E_v) \neq 0$ for all $v$.
Fixing such an $E$, and write $r_v:=\deg E_v$, so that $\sum_v r_v=r$. For 
each $e \in E(\bar{\Gamma})$ and $v$ adjacent to $e$, let $S_{(e,v)}$ be the
set of vertices in the same connected component of 
$\bar{\Gamma}\smallsetminus \{e\}$ as $v$ is, and set
$r_{(e,v)}:=\sum_{v' \in S_{(e,v)}} r_{v'}$. Thus, if $v,v'$ are the
vertices adjacent to $e$, we have $r_{(e,v)}+r_{(e,v')}=r$. We also see
that if we fix $v$, and let $e_1,\dots,e_m \in E(\bar{\Gamma})$ be the
edges adjacent to $v$, and $v_1,\dots,v_m$ the other vertices adjacent
to $e_1,\dots,e_m$, then we have
\begin{equation}\label{eq:codim-sum}
\sum_{i=1}^m r_{(e_i,v_i)} = \sum_{v' \neq v} r_{v'}=r-r_v.
\end{equation}

Now, for each $e \in E(\bar{\Gamma})$, let $v,v'$ be the adjacent 
vertices, and set $t_{(e,v)}=i$, where 
$a^{(e,v)}_{r-r_{(e,v)}}=\deg D^{(e,v)}_i$
and $i$ is critical for $D^{(e,v)}_{\bullet}$. Then set 
$t_{(e,v')}=b_{v,v'}-t_{(e,v)}$. By the multivanishing inequality, and
because $r_{(e,v)}=r-r_{(e,v')}$, 
we have $a^{(e,v')}_{r-r_{(e,v')}} \geq \deg D^{(e,v')}_{t_{(e,v')}}$.
Then the collection of all $t_{(e,v)}$ determine a multidegree 
$w \in V(\bar{G}(w_0))$ as follows: starting from some $w_{v_0}$, for
each $e$ adjacent to $v_0$, twist $t_{(e,v_0)}$ times at $(e,v_0)$.
Then traverse $\Gamma$ outward from $v_0$; for each vertex $v$ adjacent to
a previous one, and every edge $e$ not in the direction of $v_0$, twist
$t_{(e,v)}$ times at $(e,v)$. One checks that the result is a multidegree
$w$ with the property that for every $v$, and every $e$ adjacent to $v$,
to get from $w_v$ to $w$, the
number of twists required at $(e,v)$ is equal to $t_{(e,v)}$.
Then by construction we have that for each
$e$ and adjacent $v$, imposing vanishing along $D^{(e,v)}_{t_{(e,v)}}$ on 
$V^v$ results in codimension at most $r-r_{(e,v)}=r_{(e,v')}$. Thus,
the total codimension on $V^v$ is, by \eqref{eq:codim-sum}, bounded by
$r-r_v$, and $\dim V^v(-D_{w,v}) \geq r_v+1$. Thus, 
$V^v(-D_{w,v}-E_v) \neq 0$, as desired.
\end{proof}

\section{Examples and further discussion}\label{sec:examples}

We conclude by considering some examples, and discussing the question
of when one expects Brill-Noether generality for complexes of curves
or for metric graphs. The philosophy that emerges is that if we have
a given $(X_0,\bn)$, with associated complex $\fC(X_0,\bn)$ and metric
graph $\GamAB$, the complex $\fC(X_0,\bn)$ should only be Brill-Noether 
general when gluing conditions on $X_0$ automatically impose the expected 
codimension at the level of line bundle gluings (as in the case of 
pseudocompact-type curves with few nodes 
connecting pairs of components, studied in \S 5 of \cite{os23}). The
metric graph $\GamAB$ should only be Brill-Noether general when in 
addition no generality or characteristic conditions are necessary on the 
components of $X_0$ in order for $(X_0,\bn)$ to be Brill-Noether general 
with respect to limit linear series. 

Our examples have genus-$0$ components both in order to make them more
tractable and to make the comparison to metric graphs more relevant.
However, we expect the higher-genus case to behave similarly if we 
restrict our attention to the comparison between our limit linear series
and Amini-Baker limit linear series.

We first recall the definitions of tropical Brill-Noether loci.

\begin{defn}\label{def:trop-BN-2} Let $\GamAB$ be a metric graph.
Let $\Pic^d(\GamAB)$ be the set of divisors of degree $d$ up to linear
equivalence, and $W^r_d(\GamAB) \subseteq \Pic^d(\GamAB)$ the subset of
divisor classes of rank at least $r$.
\end{defn}

Then $\Pic^d(\GamAB)$ is a real torus of dimension equal to the genus
of $\GamAB$, and $W^r_d(\GamAB)$ is piecewise polyhedral in a suitable
sense; see \S 2.2, Proposition 3.7 of \cite{l-p-p1}. Thus, we define:

\begin{defn} If $\GamAB$ has genus $g$, we say it is 
\textbf{Brill-Noether general} if for all $r,d$ with $d \leq 2g-2$ we have
that $W^r_d(\GamAB)$ is empty if $\rho<0$, and has dimension
$\rho$ if $\rho \geq 0$.
\end{defn}

It is known by tropicalizing the classical case that if $\rho \geq 0$, then 
$W^r_d(\GamAB)$
always has at least some component of dimension at least $\rho$; see 
Proposition 5.1 of \cite{pf2}. On the other hand, the following 
definition, developed 
by Lim, Payne and Potashnik in \cite{l-p-p1}, is also useful
as a substitute for $\dim W^r_d(\GamAB)$ in the study of tropical 
Brill-Noether loci:

\begin{defn} For a metric graph $\GamAB$, and $r,d$, the 
\textbf{Brill-Noether rank} $w^r_d(\GamAB)$ is the largest integer $\rho'$
such that for every effective divisor $E$ on $\GamAB$ of degree $r+\rho'$,
there exists some $D$ of rank at least $r$ such that $D-E$ is effective.
\end{defn}

The following statements constitute Theorem 1.7 and Corollary 1.8 of
\cite{l-p-p1}.

\begin{thm}\label{thm:lpp} Let $R$ be a discrete valuation ring, and $K$
its field of fractions. If $X$ is a smooth projective curve over $K$
having a regular semistable model over $R$ with dual graph of the special
fiber equal to $\GamAB$, then
$$\dim W^r_d(X) \leq w^r_d(\GamAB).$$

It follows that if $\GamAB$ is an arbitrary metric graph of genus $g$, then
$$w^r_d(\GamAB) \geq \min\{\rho,g\}.$$
\end{thm}

For the sake of brevity, we will sometimes say that a divisor on a metric
graph is a ``tropical $\fg^r_d$'' or a ``$\fg^r_d$'' if it has rank $r$
and degree $d$.

We now consider our examples.

\begin{ex}\label{ex:flower} Consider the case that $X_0$ is obtained from 
rational components $Z_0, Z_1,\dots,Z_g$, with each $Z_i$ for $i>0$
glued to $Z_0$ in two nodes, and $g>2$.
\begin{figure}
\centering
\includegraphics{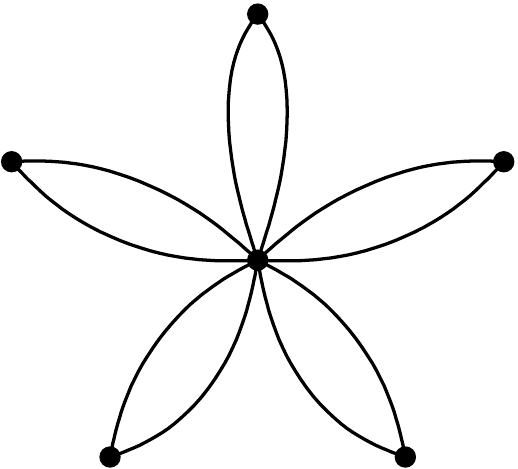}
\caption{The dual graph of Example \ref{ex:flower} in the case $g=5$.}
\label{fig:flower}
\end{figure}
Although this looks like a variation of
the compact-type curve frequently considered by Eisenbud and Harris,
with $g$ elliptic tails glued to a single rational component, the stable
model of $X_0$ is an irreducible rational curve of genus $g$, so its
behavior is rather different. In particular, generality of the points
chosen for gluing is important to the Brill-Noether generality of $X_0$. 
Corollary 5.1 of \cite{os23} implies that
in characteristic $0$, if the gluing points on $Z_0$ are general,
and the chain structure is sufficiently general,
then $X_0$ is Brill-Noether general with
respect to spaces of limit linear series on it.

In this example, we consider the narrower question of when $X_0$ supports
a limit $\fg^1_2$, and show that neither the characteristic hypothesis 
nor generality of the chain structure is relevant.
The key case to consider is that the multidegree concentrated on $Z_0$
has degree $2$ on $Z_0$, and degree $0$ elsewhere. In this case,
the $Z_0$ aspect $V^0$ of a limit $\fg^1_2$ must have multivanishing sequence
$0,2$ at each $Z_0 \cap Z_i$.
Then we see that we must have that the divisors $Z_0 \cap Z_i$ for all $i$ 
lie in a single pencil on $Z_0$,
 which is not the case if the $Z_0 \cap Z_i$ are general, since we
are assuming $g > 2$. Geometrically, the imposed condition is that if we
imbed $X_0$ as a plane conic, the lines through the images of each 
$Z_0 \cap Z_i$ must have a common intersection point.

The only other possible case is that the multidegree concentrated on $Z_0$ 
has degree
$1$ on $Z_0$, and degree $1$ on some other component, either some $Z_{i_0}$,
or some exceptional component between $Z_0$ and some $Z_{i_0}$. In this
case, we see that this multidegree is in fact concentrated at every $Z_i$
for $i \neq i_0$,
but it has degree $0$ on all except $Z_0$, so it cannot
support a pencil on any of the other components. We thus see that we 
cannot have a limit $\fg^1_2$ with such a multidegree, regardless of any
generality hypotheses.
 
Now, since the above argument never used gluing conditions, it follows from
Theorem 3.9 of \cite{he2} that we also have that $\fC (X_0,\bn)$ is 
non-hyperelliptic from the Amini-Baker perspective. More precisely,
it follows that there is no Amini-Baker limit $\fg^1_2$ on $\fC(X_0,\bn)$
which can come from a divisor supported at integral points of $\GamAB$. 
On the other hand,
since for special configurations of gluing points we do have a limit
$\fg^1_2$, we see that the associated metric graph must admit a divisor
of degree $2$ and rank $1$, regardless of genericity of edge lengths. 
Explicitly, this
is obtained simply by placing degree $2$ at the point corresponding to 
$Z_0$, and degree $0$ elsewhere. 

Finally, we address the possibility of Amini-Baker limit $\fg^1_2$s on
$\fC(X_0,\bn)$ which do not have integral support, by observing that on
$\GamAB$, any divisor of degree $2$ and rank $1$ must be linearly 
equivalent to the one described above, which has integral support. Indeed, 
we see more generally that if a given loop of $\GamAB$ doesn't have points 
on it, then to move points onto it while preserving effectivity
requires taking two from the point corresponding to $Z_0$, so specializing
to degree $2$ gives the desired assertion. 
Thus, we have that any Amini-Baker limit $\fg^1_2$ on $\fC(X_0,\bn)$ would 
have to induce the given divisor class on $\GamAB$, and is in particular 
integral. Together with the earlier argument, we can then conclude that 
$\fC(X_0,\bn)$ cannot support any Amini-Baker limit $\fg^1_2$.
\end{ex}

\begin{ex}\label{ex:binary} Consider the case that $X_0$ consists of
rational components $Z_1,Z_2$ glued to one another at $g+1$ nodes.
\begin{figure}
\centering
\includegraphics{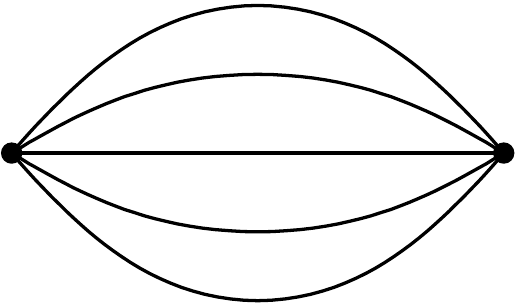}
\caption{The dual graph of Example \ref{ex:binary} in the case $g=4$.}
\label{fig:binary}
\end{figure}
This situation is studied in \S 7 of \cite{os23}, where in particular
Corollary 7.4 shows that if the gluing points are sufficiently general,
then $X_0$ has maximal gonality with respect to our notion of limit linear
series (independent of choice of chain structure or enriched structure).
On the other hand, if the gluing points are too special, then $X_0$ has
a limit $\fg^1_2$, and importantly, the failure of generality affects not
the behavior of individual components, but the transversality of the gluing
condition. Consequently, it is not surprising that even for general
gluing points, $\fC(X_0,\bn)$ admits an Amini-Baker limit $\fg^1_2$,
namely obtained by placing one point each (anywhere) on $Z_1$ and $Z_2$,
and taking spaces of rational functions to move the chosen point to any
other point.

Of course, it then follows that the associated metric graph also always
supports a tropical $\fg^1_2$, in this case obtained by placing degree
$1$ on each of the points corresponding to $Z_1$ and $Z_2$.
\end{ex}

In \S 5 of \cite{os23}, we studied curves of pseudocompact type with
at most three nodes connecting any given pair of components, and suitably
general chain structures. We showed that in this case, gluing conditions
automatically impose the expected codimension, and we observed that in the
case of a chain of loops, the condition we recovered on chain structures
precisely matched the genericity condition of Cools, Draisma, Payne and
Robeva in \cite{c-d-p-r}. This motivated us to ask in Question 5.4 
whether chains with at most three edges connecting adjacent nodes (and
generic edge lengths) are Brill-Noether general. We also asked in 
Question 5.5 whether it is possible for any other metric multitrees to
be Brill-Noether general. The above examples suggest that the latter 
question ought to have a negative answer. 

This perspective is amplified by work of Kailasa, Kuperberg and Wawrykow 
\cite{k-k-w1} showing that the only trees of loops which can be Brill-Noether
general are chains of loops. They also show that a tree of loops fails to
have maximal gonality unless it is obtained from a chain of genus one 
less by addition of a single loop. We use the techniques of our examples 
to prove gonality bounds complementing those of \cite{k-k-w1}. Following
their lead, we also study dimensions of Brill-Noether loci in the remaining 
cases, proving in particular in Corollary \ref{cor:bad-multitrees}
that the only metric multitrees which can be Brill-Noether general are the 
above-mentioned chains with at most three edges between adjacent vertices.
Specifically, the following two propositions generalize the techniques of 
our examples to study disconnecting (multi)edges and disconnecting vertexes
of metric graphs. They say roughly that if a metric graph has a vertex 
which disconnects it into three or more components, or two adjacent 
vertices which are connected by four or edges, and disconnect the
graph, then we do not have maximal gonality. More precisely, we show that
we can only have maximal gonality in very specific circumstances, and that
in these cases, the dimension of $W^1_{\left\lceil\frac{g}{2}\right\rceil+1}$ 
is nonetheless larger than expected.

Because separating edges (including ``edges to nowhere'') have no effect on 
the Brill-Noether theory of a graph,
to simplify the following statement we assume that $\GamAB$ has none.

\begin{prop}\label{prop:low-gonality-1} Let $\GamAB$ be a metric graph
without separating edges, and
$v \in \GamAB$ a point which disconnects $\GamAB$. Let 
$\GamAB_1,\dots,\GamAB_n$ be the closures of the connected components of
$\GamAB \smallsetminus \{v\}$, so that $\GamAB$ is the wedge of the
$\GamAB_i$ at $v$. Let 
$n_1$ the number of $\GamAB_i$ with odd genus, and $n_2$ the number with 
even genus, and suppose that $n > 2$. Then $\GamAB$ does not
have maximal gonality, unless the following conditions are all satisfied:
$n=3$, $n_1=2$, $n_2=1$, at least
one of the $\GamAB_i$ has genus precisely $1$, and the $\GamAB_i$ with
even genus has maximal gonality.

Furthermore, even when the above conditions are satisfied, if $g$ denotes 
the genus of
$\GamAB$, the dimension of $W^1_{\left\lceil\frac{g}{2}\right\rceil+1}(\GamAB)$
is positive.
\end{prop}

Note that in the last statement, 
$W^1_{\left\lceil\frac{g}{2}\right\rceil+1}(\GamAB)$ has expected dimension 
$0$.

\begin{proof}
Set $\epsilon=1$ if $g$ is odd, and $\epsilon=0$
if $g$ is even. Then
\begin{multline*}
\left\lceil \frac{g}{2}\right\rceil + 1 = \frac{g+\epsilon}{2}+1
=1+\frac{\epsilon}{2} + \sum_i \frac{g_i}{2} 
=1+\frac{\epsilon-n_1}{2} + \sum_i \left\lceil\frac{g_i}{2}\right\rceil
\\=1-n+\frac{\epsilon-n_1}{2} 
+ \sum_i (\left\lceil\frac{g_i}{2}\right\rceil+1).\end{multline*}
If we take a pencil of degree 
$\left\lceil\frac{g_i}{2}\right\rceil+1$ on each
$\GamAB_i$, we may assume by Theorem \ref{thm:lpp} it contains $v$ with 
multiplicity at least $1$
when $g_i$ is even, and with multiplicity at least $2$ when $g_i$ is odd.
Thus, if we take the (additive) ``least common multiple''
of these divisors on $\GamAB$, we obtain a pencil of degree at most
$$2+\left(\sum_i (\left\lceil\frac{g_i}{2}\right\rceil+1)\right)-n-n_1
=\left\lceil \frac{g}{2}\right\rceil+1+(1-\frac{\epsilon+n_1}{2}).$$
Since $\epsilon=1$ if and only if $n_1$ is odd, 
$\frac{\epsilon+n_1}{2}=\left\lceil\frac{n_1}{2}\right\rceil$, so we 
conclude that the above is strictly smaller than
$\left\lceil \frac{g}{2}\right\rceil + 1$ precisely when $n_1$ is at least 
$3$. Thus, $\GamAB$ does not have maximum gonality when $n_1 \geq 3$. 
Note also that if any even-genus $\GamAB_i$ fails to have maximal gonality,
we can carry out the above construction with one smaller total degree,
so that in this case we find that $\GamAB$ doesn't have maximal gonality
as long as $n_1 \geq 1$.

On the other hand, if we take pencils of degree 
$\left\lceil\frac{g_i}{2}\right\rceil+2$ on each $\GamAB_i$ with $g_i$ even, 
we may assume that these contain $v$ with multiplicity at least $3$. Then
taking the lcm as above, we obtain degree at most
$$3+\left(\sum_i (\left\lceil\frac{g_i}{2}\right\rceil+1)\right)+n_2-2n-n_2
=\left\lceil \frac{g}{2}\right\rceil+1+(2-\frac{\epsilon+n_1}{2}-n_2).$$
As above, $\GamAB$ does not have maximum gonality when 
$\left\lceil \frac{n_1}{2}\right\rceil+n_2>2$. 
In particular, we can have maximum gonality only when $n_2 \leq 2$.
Putting these together, we see that if we assume $n \geq 3$, we can have 
maximum gonality only when $n_1 =2$ and $n_2 = 1$. Together with our
earlier observation, we also see that we need maximal gonality on the
$\GamAB_i$ with (positive) even genus.

Furthermore, suppose that $g_i>1$ for all $i$, so that each odd $g_i$ is 
at least $3$. We can then take pencils of degree
$\left\lceil\frac{g_i}{2}\right\rceil+2$ on each $\GamAB_i$ with $g_i$ odd, 
and choose divisors in these pencils with multiplicity at least $4$ at $v$.
Taking the lcm as above, we obtain a pencil of degree at most 
$$4+\left(\sum_i (\left\lceil\frac{g_i}{2}\right\rceil+1)\right)
+n_2+n_1-3n-n_1
=\left\lceil \frac{g}{2}\right\rceil+1+(3-\frac{\epsilon+ n_1}{2}-n),$$
so the resulting gonality is less than maximal if 
$\left\lceil \frac{n_1}{2}\right\rceil+n>3$, which it is in the case $n_1=2$ 
and $n_2=1$. We thus conclude the first assertion of the proposition.

Now, suppose that we have $n_1=2$ and $n_2=1$, and suppose without loss
of generality that $g_1$ and $g_2$ are odd, and $g_3$ is even. Then consider 
the construction in the second paragraph above, with pencils of degree
$\left\lceil\frac{g_i}{2}\right\rceil+1$ on $\GamAB_i$ for $i=1,2$,
and $\left\lceil\frac{g_3}{2}\right\rceil+2$ on $\GamAB_3$. We see
that if instead of requiring the latter to have degree $3$ at $v$, we
require it to have degree $2$, we are free to also impose degree $1$ at
any other point $v'$ on $\GamAB_3$. In this case, the lcm has degree
$$2+ \left(\sum_i (\left\lceil\frac{g_i}{2}\right\rceil+1)\right)+1-6
\left\lceil \frac{g}{2}\right\rceil+1,$$
so for each choice of $v'$, we obtain a pencil of degree
$\left\lceil \frac{g}{2}\right\rceil+1$ on $\GamAB$. 
We claim that the constructed pencils cannot all be linearly equivalent: 
indeed, if not their restrictions to $\GamAB_3$
would likewise have to be linearly equivalent,
which would imply, factoring out the base points
at $v$, that $\GamAB_3$ has a pencil of degree 
$\left\lceil\frac{g_3}{2}\right\rceil$.
This contradicts our maximal gonality hypothesis.
Thus, we have that 
$W^1_{\left\lceil\frac{g}{2}\right\rceil+1}(\GamAB)$ is 
positive-dimensional, as claimed.
\end{proof}

\begin{prop}\label{prop:low-gonality-2} Let $\GamAB$ be a metric graph, and
suppose that there are points $v_1,v_2 \in \GamAB$ and $\GamAB_1,\GamAB_2
\subseteq \GamAB$ disjoint and containing $v_1,v_2$ respectively such that 
$\GamAB$ is obtained from $\GamAB_1 \cup \GamAB_2$ by adding $m$ disjoint
paths from $v_1$ to $v_2$. Suppose further that $m \geq 4$. Then
$\GamAB$ does not have maximal gonality, unless both $\GamAB_i$ have maximal 
gonality, and we either have $m=4$ with at least one $\GamAB_i$ having odd 
genus, or $m=5$ with both $\GamAB_i$ having odd genus.

Furthermore, even when the above conditions are satisfied, if $g$ denotes 
the genus of
$\GamAB$, the dimension of 
$W^1_{\left\lceil\frac{g}{2}\right\rceil+1}(\GamAB)$
is strictly greater than the expected 
$2\left\lceil\frac{g}{2}\right\rceil -g$.
\end{prop}

\begin{proof}
Let $g_1,g_2$ be the genera of $\GamAB_1,\GamAB_2$, and for
$i=1,2$, let $D_i$ be a $\fg^1_{\left\lceil\frac{g_i}{2}\right\rceil+1}$ on $\GamAB_i$
containing $v_i$ in its support. Then it follows that $D_1+D_2$ is a
$\fg^1_{\left\lceil\frac{g_1}{2}\right\rceil+\left\lceil\frac{g_2}{2}\right\rceil+2}$ on
$\GamAB$. We have that $\GamAB$ has genus $g:=g_1+g_2+m-1$, so if neither
$g_i$ is odd, or if only one is odd and $m = 5$, or if $m \geq 6$, we have that
$$\left\lceil\frac{g}{2}\right\rceil+1 > 
\left\lceil\frac{g_1}{2}\right\rceil+\left\lceil\frac{g_2}{2}\right\rceil+2,$$ 
as desired. Moreover, even if $m=4$ and one of the $g_i$ is odd or $m=5$
and both of the $g_i$ are odd, we see that we have
$$\left\lceil\frac{g}{2}\right\rceil+1 = 
\left\lceil\frac{g_1}{2}\right\rceil+\left\lceil\frac{g_2}{2}\right\rceil+2,$$ 
so if either of the $\GamAB_i$ have less than maximal gonality, we can
still reduce the degree in the above construction strictly below 
$\left\lceil\frac{g}{2}\right\rceil+1$.

Next, suppose that both $\GamAB_i$ have maximal gonality, and either 
$m=4$ and at least one of the $g_i$ is odd, or $m=5$ and both
$g_i$ are odd. First suppose that $\GamAB$ has even genus, so
that either $m=4$ and exactly one $g_i$ is odd, or $m=5$
and both $g_i$ are odd. We claim that
$W^1_{\left\lceil \frac{g}{2}\right\rceil+1}$ is infinite. Without loss of generality,
in either case suppose that $g_2$ is odd. Fix $D_1$ a 
$\fg^1_{\left\lceil\frac{g_1}{2}\right\rceil+1}$ on $\GamAB_1$ having degree at least
$1$ at $v_1$. Next, for any $v' \in \GamAB_2$, there exists $D_{2,v'}$ a 
$\fg^1_{\left\lceil\frac{g_2}{2}\right\rceil+1}$ on $\GamAB_2$ containing $v_2+v'$.
These cannot all be linearly dependent on $\GamAB_2$, as otherwise we 
could factor out $v_2$ and obtain a
$\fg^1_{\left\lceil\frac{g_2}{2}\right\rceil}$ on $\GamAB_2$.
This implies that the $D_1+D_{2,v'}$ are not all linearly equivalent on
$\GamAB$, 
so $W^1_{\left\lceil \frac{g}{2}\right\rceil+1}$ is infinite, as desired.

Finally, we want to show that if $m=4$ and both $g_i$ are odd, we have that
$W^1_{\left\lceil \frac{g}{2}\right\rceil+1}$ is at least $2$-dimensional. As above,
construct families $D_{i,v'_i}$ of 
$\fg^1_{\left\lceil\frac{g_i}{2}\right\rceil+1}$s on the $\GamAB_i$ containing 
$v_i+v'_i$, where $v'_i$ varies on $\GamAB_i$. By the maximal gonality
condition, we see that we cannot have all the constructed pencils on either 
$\GamAB_i$ being linearly equivalent to one another. But now we observe
that given $d_1,d_2$, if we set $d=d_1+d_2$, the map 
$$\Pic^{d_1}(\GamAB_1) \times \Pic^{d_2}(\GamAB_2) \to \Pic^d(\GamAB)$$
given by adding divisors is injective: indeed, any piecewise linear
function giving a linear equivalence on $\GamAB$ has to be constant on
the paths connecting $v_1$ to $v_2$, since degree is preserved on
$\GamAB_1$ and $\GamAB_2$, and hence it induces linear equivalences on
$\GamAB_1$ and $\GamAB_2$. Moreover, this map is expressed as the
composition of ``linear'' maps $\Pic^{d_i}(\GamAB_i) \to \Pic^{d_i}(\GamAB)$
(in the sense that these are induced by linear projections on the universal 
covering spaces) composed with the addition map studied in Lemma 3.5 of
\cite{l-p-p1}, both of which behave well with respect to polyhedral 
structures. It thus follows that
$W^1_{\left\lceil \frac{g}{2}\right\rceil+1}$ is at least $2$-dimensional, as
desired.
\end{proof}

Putting together Propositions \ref{prop:low-gonality-1} and
\ref{prop:low-gonality-2}, we obtain the following expected negative 
result.

\begin{cor}\label{cor:bad-multitrees} A metric multitree can only be
Brill-Noether general if it is a multichain, with at most three edges
between adjacent vertices.
\end{cor}

Beyond Corollary \ref{cor:bad-multitrees}, there are other factors which
suggest that Brill-Noether generality should be relatively rare for
graphs. For instance, even in the compact-type case it is not true that
spaces of limit linear series are well behaved in general in positive
characteristic, even when the components of the curves are general.
Since Brill-Noether generality of graphs is characteristic-independent,
this places further constraints on the possibilities for Brill-Noether
generality. Taken together with our examples and propositions, this 
may suggest that very few graphs should be Brill-Noether general. However,
once one moves beyond multitrees, it is unclear what to expect. On the one
hand, if a graph is heavily interconnected it is harder to concentrate
degrees on particular vertices, which for instance may mitigate difficulties
arising from positive characteristic. On the other hand, Jensen has given
an example \cite{je1} of a trivalent graph without any separating edges 
which fails to be Brill-Noether
general regardless of edge length. Ultimately, further development of our
understanding of limit linear series beyond the pseudocompact type case
should also help to guide intuition on which graphs ought to be
Brill-Noether general.

\appendix

\section{Background on the Amini-Baker theory}

In this appendix, we describe background results on Amini-Baker limit
linear series. We begin with the fact that restricted rank of a divisor
on a metrized complex can be checked using only points of the curves
$C_v$. In Theorem A.1 of \cite{a-b1}, Amini and Baker prove a stronger
statement in the case of non-restricted rank, but this statement fails
for restricted rank: indeed, already in the case of incomplete linear
series on smooth curves, no fixed finite set of points can suffice to
check the rank of every linear series (see also Example \ref{ex:frob} 
below). 

\begin{thm}\label{thm:a-b-r-rk} Let $\fC$ be a metrized complex, $\cD$ a 
divisor on $\fC$, and $(H_v)_v$ a tuple of $(r+1)$-dimensional subspaces
of the $K(C_v)$. Fix $\cR$ any collection of points of $\bigcup_v C_v$
containing infinitely many points from each $C_v$. Then the restricted 
rank of $\cD$ is at least $r'$
if and only if for every effective divisor $\cE$ on $\fC$ of degree $r'$
supported only on $\cR$, there exists a rational function $\ff$ on
$\fC$ with $f_v \in H_v$ for all $v$, and such that $\cD-\cE+\dv  \ff$
is effective.
\end{thm}

Recall that we have a running hypothesis that the graph $\GamAB$ 
underlying $\fC$ is loopless.

We will carry out a close analysis of linear equivalence of divisors on
metrized complexes for divisors satisfying the analogue of our 
admissible multidegree condition:

\begin{defn} We say a divisor $\cD$ on a metrized complex $\fC$ is
\textbf{edge-reduced} if it is effective on $\GamAB\smallsetminus V(\GAB)$,
and if further $\cD$ has degree at most $1$ on every connected component
of $\GamAB\smallsetminus V(\GAB)$.
\end{defn}

Now, note that given any edge-reduced
divisor $\cD$ and any values $(c_v)_{v \in V(\GAB)}$, 
there is a unique
continuous piecewise linear function on $\GamAB$, which we denote by 
$f_{\cD,c_{\bullet}}$, such that $f_{\cD,c_{\bullet}}(v)=c_v$ for all $v$,
and such that $\cD+\dv f_{\cD,c_{\bullet}}$ is still edge-reduced 
(in particular, if $\cD-\cD'=\dv \ff$, and $\cD'$ is also edge-reduced,
then $f_{\GamAB}$ is necessarily of this form). If an edge of $\GamAB$ 
connects $v$ to $v'$, and $c_v>c_{v'}$, this has the effect of ``moving a 
chip distance $f(v)-f(v')$'' from $C_v$ toward $C_{v'}$ along the edge in 
question (here, each time the chip reaches $C_{v'}$, a new chip must be
moved starting from $C_v$, possibly creating a negative coefficient at
the corresponding point).
Note that this construction is additive in $c_{\bullet}$, in the sense
that if we have $\cD$ and $c_{\bullet},c'_{\bullet}$, then
$$f_{\cD,c_{\bullet}+c'_{\bullet}}=f_{\cD,c_{\bullet}}+f_{\cD', c'_{\bullet}},$$
where $\cD'=\cD+\dv f_{\cD,c_{\bullet}}$. Also note that if we are given
$\cD$ edge-reduced and $c_{\bullet}$ such that $\cD$ and 
$\cD+\dv f_{\cD,c_{\bullet}}$ are both effective, then
$\cD+\dv f_{\cD,\alpha c_{\bullet}}$ will also be effective 
(and edge-reduced) for all $\alpha \in [0,1]$.

\begin{lem}\label{lem:equiv-decompose} Let $\fC$ be a metrized complex,
and suppose that $\cD$ and $\cD'$ are linearly equivalent edge-reduced
effective divisors on $\fC$, with $\cD-\cD'=\dv \ff$. Then there exists
an ordering $v_1,\dots,v_n$ on the vertices of $\fC$ and tuples
$c^i_{\bullet}$ for $i=1,\dots,n-1$ such that if we set 
$\cD_1=\cD+\dv f_{v_1}$, and for $i=2,\dots,n$ write
$$\cD_i=\cD_{i-1}+\dv f_{\cD_{i-1},c^{i-1}_{\bullet}}+\dv f_{v_i},$$
the following conditions are satisfied:
\begin{ilist}
\itm $\cD_{i}+\dv f_{\cD_{i},\alpha c^{i}_{\bullet}}$ is effective and
edge-reduced for $i=1,\dots,n-1$ and $\alpha \in [0,1]$;
\itm if for some $i_1<i_2$ we have $f_{\GamAB}(v_{i_1})=f_{\GamAB}(v_{i_2})$,
then $f_{\cD_{i},c^i_{\bullet}}$ is constant for $i_1 \leq i < i_2$;
\itm for each $i=1,\dots,n$, we have $\cD_i-\dv f_{v_{i}}-\cD$ effective 
on $C_{v_j}$ for all $j \geq i$, and $\cD_i-\cD'$ effective on $C_{v_j}$ 
for all $j \leq i$;
\itm $\cD_n=\cD'$.
\end{ilist}
\end{lem}

\begin{proof} 
Choose the ordering $v_1,\dots,v_n$ so that
$f_{\GamAB}(v_i) \leq f_{\GamAB}(v_{i-1})$ for $i=2,\dots,n$. Note that
being edge-reduced implies by hypothesis effectivity away from the 
$C_{v_i}$, so since our constructed divisors will always be edge-reduced,
in order to check effectivity it suffices to consider 
the $C_{v_i}$ one at a time. For $j=1,\dots,n-1$, we will set
$c^j_i=f_{\GamAB}(v_{j})$ for $i \leq j$, and 
$c^j_i=f_{\GamAB}(v_{j+1})$ for $i > j$.
Then observe that for each $i=1,\dots,n$, we have 
$$\sum_{j=1}^{n-1} c^j_i 
= f_{\GamAB} (v_i) + \sum_{\ell=2}^{n-2} f_{\GamAB} (v_{\ell}),$$
so we conclude that $\cD_n=\cD'$ by the additivity of the 
$f_{\cD,c_{\bullet}}$ construction.

Now, $\dv f_{\GamAB}$ is anti-effective on $C_{v_1}$, so effectivity of
$\cD'$ implies effectivity of $\cD_1$ of $C_{v_1}$, and effectivity of
$\cD$ implies that $\cD_1$ is effective on the rest of $\fC$. 
Next, $\dv \left(f_{\GamAB}-f_{\cD_1,c^1_{\bullet}}\right)$ is
anti-effective on both $C_{v_1}$ and $C_{v_2}$, and 
$\dv f_{\cD_1,c^1_{\bullet}}$ is effective on $C_{v_i}$ for $i>1$, so as 
before the effectivity 
of $\cD_1$ and $\cD'$ implies first that $\cD_1+\dv f_{\cD_1,c^1_{\bullet}}$ 
is effective, and second that $\cD_2$ is effective.
From the former, we also conclude that
$\cD_1+\dv f_{\cD_1,\alpha c^1_{\bullet}}$ is effective for all
$\alpha \in [0,1]$. The same argument works for all $j$:
we have that $\dv f_{\GamAB}-\sum_{i=1}^j \dv f_{\cD_i,c^i_{\bullet}}$ is 
anti-effective on $C_{v_1},\dots,C_{v_{j+1}}$, and 
$\dv f_{\cD_j,c^j_{\bullet}}$ is effective on $C_{v_{j+1}},\dots,C_{v_n}$,
so the effectivity of $\cD_j$ and $\cD'$ implies that
$\cD_j +\dv f_{\cD_j,c^j_{\bullet}}$ and $\cD_{j+1}$ are both effective.
We then conclude that assertion (i) is satisfied. Assertion (ii) is clear
from the construction. Finally, assertion (iii) follows from the facts that 
$\sum_{\ell=1}^{i-1} \dv f_{\cD_{\ell},c^{\ell}_{\bullet}}$ is effective on 
$C_{v_{j}}$ for all $j\geq i$, and that 
$\sum_{\ell=i}^{n-1} \dv f_{\cD_{\ell},c^{\ell}_{\bullet}}$ 
is anti-effective on $C_{v_j}$ for all $j \leq i$.
\end{proof}

\begin{proof}[Proof of Theorem \ref{thm:a-b-r-rk}]
First, the inductive portion of the
proof of Lemma A.3 of \cite{a-b1} goes through verbatim in the case
of restricted rank, showing that if we fix the $(H_v)_v$, and a subset
$\cR$ of $\bigcup_v C_v$, then $\cR$ suffices to test whether the restricted 
rank of $\cD$ is at least $1$ for every $\cD$ if and only if $\cR$ in fact 
suffices to test the restricted rank of every $\cD$.

Now, given $\cD_0$, suppose that $\cD_0-Q$ is $(H_v)_v$-equivalent to an
effective divisor for all $Q \in \cR$. We first claim that in fact we must
have that $\cD_0-Q$ is $(H_v)_v$-equivalent to an effective divisor for all
$Q \in \bigcup_v C_v$. For each $Q \in \cR$, let $\ff_Q$ be a rational
function such that $\cD_0-Q+\dv \ff_Q$ is effective, and observe that as
$Q$ varies, for any given $v$ there are only finitely many possibilities
for the $C_v$-part of $\dv((f_Q)_{\GamAB})$.
In particular, for any given $v$, there exists an infinite set $\cR'_v$ of 
points $Q \in C_v \cap \cR$ such that all the $C_v$-parts of the 
$\dv((f_Q)_{\GamAB})$
agree; denote these common divisors on $C_v$ by $\cD'_v$. Thus we
have $(\cD_0)_v-Q+\cD'_v+\dv (\ff_Q)_v$ effective on $C_v$ for all 
$Q \in \cR'_v$, and we may also without loss of generality assume that
the $\ff_Q$ for $Q \in \cR'_v$ all agree of $\GamAB$ and on $C_{v'}$ with
$v' \neq v$; denote these fixed rational functions by $f_{\GamAB}$ and
$f_{v'}$ for $v' \neq v$.
Since $\cR'_v$ is infinite, we conclude that $(\cD_0)_v+\cD'_v+\dv(H_v)$
must contain a pencil, so that for all $Q \in C_v$, we have some
$f_{v,Q} \in H_v$ with $(\cD_0)_v-Q+\cD'_v+\dv (\ff_Q)_v$ effective on $C_v$.
Extending $f_{v,Q}$ to a rational function on $\fC$ using our fixed
$f_{\GamAB}$ and $f_{v'}$ for $v' \neq v$, we see $\cD_0-Q$ is
$(H_v)$-equivalent to an effective divisor, as claimed.

Next, suppose we are given
$P \in \GamAB \smallsetminus V(\GAB)$. Let $v,v'$ be the endpoints of the
edge $e$ of $\GamAB$ containing $P$, and $Q,Q'$ the corresponding points
in $\fC$. 
By hypothesis, there exist effective divisors $\cD,\cD'$
which are $(H_v)_v$-equivalent to $\cD_0$ and which contain $Q$ and $Q'$
in their support, respectively. Without loss of generality, we may assume
that $\cD$ and $\cD'$ are both edge-reduced. Fix $\ff$ such that
$\cD-\cD'=\dv \ff$, and let $v_1,\dots,v_n$, 
$c^i_{\bullet}$ and $\cD_i$ be as in Lemma 
\ref{lem:equiv-decompose}. Suppose $v=v_{i_1}$ and $v'=v_{i_2}$.
Switching the roles of $Q$ and $Q'$ if necessary, we may assume that
$i_1<i_2$.

We claim that (at least) one of the following must occur:
\begin{alist}
\itm for some $i$ and $\alpha$, the divisor
$\cD_{i-1}+\dv f_{\cD_{i-1},\alpha c^{i}_{\bullet}}$ contains $P$;
\itm $\cD$ contains a point $P'$ (not necessarily strictly) between $P$ and 
$Q'$ on $e$;
\itm $\cD'$ contains a point $P'$ (not necessarily strictly) between $Q$ and 
$P$ on $e$;
\itm $\cD_{i_1}-\dv f_{v_{i_1}}$ contains $Q$, and $\cD_{i_2}$ contains $Q'$,
and $\dv f_{\cD_{i},c^{i}_{\bullet}}=0$ for $i_1 \leq i < i_2$.
\end{alist}
In case (a), we have constructed an effective divisor 
$(H_v)_v$-equivalent to $\cD_0$ and containing $P$, while in cases (b) and (c)
we can easily modify $\cD$ or $\cD'$ by the divisor of a 
function on $\GamAB$, constant away from $e$, such that the result contains 
$P$. Finally, in case (d), note that by construction 
$\cD_{i_1}-\dv f_{v_{i_1}}$ is effective, and because
$\dv f_{\cD_{i},\alpha c^{i}_{\bullet}}=0$ for $i_1 \leq i < i_2$, we have
that $\cD_{i_2}$ and $\cD_{i_1}$ differ only by the $\dv f_{v_i}$ for
$i_1<i \leq i_2$. It follows
$\cD_{i_2}-\dv f_{v_{i_1}}$ is effective and contains both $Q$ and $Q'$,
so as above, we can modify it by the divisor of a function on $\GamAB$
to obtain an effective divisor containing $P$.

Thus, we have reduced to proving the claim. Because $i_1<i_2$, we have
$f_{\GamAB}(v) \geq f_{\GamAB}(v')$; if we write 
$\mu=f_{\GamAB}(v)-f_{\GamAB}(v')$, then the effect of $\dv f_{\GamAB}$
is to move a chip a distance of $\mu$ along $e$. Write $\ell$ for the
length of $e$, and if $\cD$ has a chip on the interior of $e$, let $p$ be 
the distance from $Q$ to that chip. Otherwise, set $p=0$. Also, let $p'$ be
the distance from $Q$ to $P$. Then 
$\dv f_{\GamAB}$ has the effect of moving a chip from distance $p$ to 
distance $p+\mu$, considered modulo $\ell$. By Lemma 
\ref{lem:equiv-decompose} (i), if $p'$ is within
this range, then we are in case (a) above. But if $p'$ is not within
this range, then we must have either $p > p'$ or $p+\mu < p'$.
The former gives us case (b) above, while the latter almost implies
that we are in case (c). The only special case to consider is
that if $p+\mu=0$ (modulo $\ell$), in which case it is not \textit{a priori}
the case that $\cD'$ has a chip at the corresponding position, which is $Q$. 
However, if we are not in case (b) we can in particular suppose that
$p < p'$, so if either $p$ or $\mu$ is positive, then $p'$ is in the
range from $p$ to $p+\mu$ modulo $\ell$, and we are in case (a). Thus,
it remains to consider the case that $p=\mu=0$. In this case
$f_{\GamAB}(v)=f_{\GamAB}(v')$, so according to Lemma 
\ref{lem:equiv-decompose} (ii) and (iii), we have that
$\dv f_{\cD_{i-1},c^i_{\bullet}}=0$ for $i_1 \leq i < i_2$, and we
also have that $\cD_{i_1}-\dv f_{v_{i_1}}-\cD$ is effective on $C_{v}$,
and $\cD_{i_2}-\cD'$ is effective on $C_{v'}$. In particular, 
$\cD_{i_1}-\dv f_{v_{i_1}}$ contains $Q$ and $\cD_{i_2}$ contains $Q'$,
so we are in case (d) above, proving the claim and the theorem.
\end{proof}

\begin{rem} We need Theorem \ref{thm:a-b-r-rk} for our analysis
of Amini-Baker limit linear series (specifically, in the proofs of 
Theorem \ref{thm:forget} and Proposition \ref{prop:forget-pseudocompact},
treating our forgetful map construction). In fact, Theorem \ref{thm:a-b-r-rk}
is also used implicitly in \cite{a-b1}, in the proofs of Theorems 5.9 and 
5.11, where they consider restricted ranks in the context of their limit 
linear series.
\end{rem}



\begin{ex}\label{ex:frob} For complete linear series, it is classical
(and an easy consequence of the Riemann-Roch theorem)
that if $X$ has genus $g$, and we fix any $g+1$ points $P_i$ on $X$, then
$r(D) \geq r'$ if and only if for every effective $E$ of degree $r'$
supported on the $P_i$, we have that $D-E$ is linearly equivalent to an
effective divisor. 
This is what is generalized in Theorem A.1 of
\cite{a-b1} to metrized complexes. This statement fails trivially for
incomplete linear series, as we could for instance take a $\fg^0_d$
with base points at all the $P_i$. More interestingly, even if we fix
a $\fg^r_d$ in advance, there may not exist any choice of $g+1$ points
with the above rank-determining property. Indeed, if we consider
the Frobenius map on $\PP^1$, it is a $\fg^1_p$. If we choose any $P$
on $\PP^1$, we will have that $pP$ is a representative of our $\fg^1_p$,
so $pP-iP$ is effective for all $i \leq p$, and if $\{P\}$ were 
rank-determining, this would contradict that our linear series had rank
$1$. On the other hand, we do see in this example that if we choose
any $P_1 \neq P_2$, since our linear series has no divisor supported
at both $P_1$ and $P_2$, the set $\{P_1,P_2\}$ is rank-determining for
the given linear series.
\end{ex}

We also give a basic statement on the specialization of Amini-Baker limit
linear series, saying in essence that specialization is compatible in
the obvious way with resolution of singularities.

\begin{prop}\label{prop:specialize-ab} Let $B$ be the spectrum of a discrete
valuation ring, and $\pi:X \to B$ flat and proper with nodal special fiber
$X_0$ having smooth components, and smooth generic fiber $X_{\eta}$. 
Let $\widetilde{\pi}:\widetilde{X} \to B$ be the regularization of 
$X$, let $\Gamma$ and $\widetilde{\Gamma}$ be the corresponding dual graphs,
and let $\fC$ and $\widetilde{\fC}$ be the induced metrized complexes 
of curves.

If $(\sL_{\eta},V_{\eta})$ is a $\fg^r_d$ on $X_{\eta}$, let
$(\cD,(H_v)_v)$ and $(\widetilde{\cD},(\widetilde{H}_v)_v)$ be the
Amini-Baker limit $\fg^r_d$s on $\fC$ and $\widetilde{\fC}$ respectively,
induced by specialization of $(\sL_{\eta},V_{\eta})$ 
(Theorem 5.9 of \cite{a-b1}). Then $H_v=\widetilde{H}_v$ for all 
$v \in V(\Gamma)$, and $\cD$ is obtained from $\widetilde{\cD}$ by,
for each $v \in V(\widetilde{\Gamma}) \smallsetminus V(\Gamma)$,
replacing $D_v$ with $\deg D_v$ times the point of $\GamAB$ corresponding
to $v$.
\end{prop}

\begin{proof} To see that $\cD$ is as described, we need that the inclusion
of $\GamAB$ into $X_{\eta}^{\an}$ and the retraction map $X^{\an} \to \GamAB$
are both unchanged under the replacement of $X$ with $\widetilde{X}$.
For the former, see Theorem 4.11 and Proposition 4.21 (2) of \cite{b-p-r2},
while the latter follows from the purely topological nature of the
retraction map (see Definition 3.7 of \cite{b-p-r2}). To see that
$H_v$ is as described, we use that the identification of the function
field of $C_v$ with that of $C_x$ (where $x$ is the relevant point
of $X_{\eta}^{\an}$) described in Remark 4.18 of \cite{b-p-r2} is also
canonical, coming as it does from viewing $x$ as a divisorial valuation
on the function field of $X_{\eta}$.
\end{proof}

\bibliographystyle{amsalpha}
\bibliography{gen}

\end{document}

%% file: headers.tex
\usepackage{amssymb,euscript,amsmath, mathrsfs, amscd, mathabx}
\usepackage[pdftex]{graphicx}
\usepackage[pdftex]{color}
\usepackage{bm}

\newcounter{ENUM}
\newcommand{\itm}{\item}
\newenvironment{ilist}[1][0]{\renewcommand{\theENUM}{\roman{ENUM}}\renewcommand{\itm}{\addtocounter{ENUM}{1}\item[(\theENUM)]}\begin{itemize}\setcounter{ENUM}{#1}}{\end{itemize}}
\newenvironment{Ilist}{\renewcommand{\theENUM}{\Roman{ENUM}}\renewcommand{\itm}{\addtocounter{ENUM}{1}\item[(\theENUM)]}\begin{itemize}\setcounter{ENUM}{0}}{\end{itemize}}
\newenvironment{alist}[1][0]{\renewcommand{\theENUM}{\alph{ENUM}}\renewcommand{\itm}{\addtocounter{ENUM}{1}\item[(\theENUM)]}\begin{itemize}\setcounter{ENUM}{#1}}{\end{itemize}}
\newenvironment{nlist}[1][0]{\renewcommand{\theENUM}{\arabic{ENUM}}\renewcommand
{\itm}{\addtocounter{ENUM}{1}\item[(\theENUM)]}\begin{itemize}\setcounter{ENUM}{
#1}}{\end{itemize}}

\newcommand{\margh}[1]{}

\def\risom{\overset{\sim}{\rightarrow}}

\input xy
\xyoption{all}
\CompileMatrices

\def\ZZ{{\mathbb Z}}

\def\PP{{\mathbb P}}

\def\bn{{\bm{n}}}
\def\cD{{\mathcal D}}
\def\cE{{\mathcal E}}

\def\cA{{\mathcal A}}

\def\cR{{\mathcal R}}

\def\sL{{\mathscr L}}

\def\sO{{\mathscr O}}

\def\fC{{\mathfrak C}}
\def\ff{{\mathfrak f}}
\def\fg{{\mathfrak g}}

\def\e{\varepsilon}
\def\vp{\varphi}

\def\dv{\operatorname{div}}

\def\Pic{\operatorname{Pic}}

\def\ord{\operatorname{ord}}

\def\can{\operatorname{can}}

\def\an{\operatorname{an}}
\def\GAB{\widehat{G}}
\def\GamAB{\widehat{\Gamma}}

\newtheorem{thm}{Theorem}[section]
\newtheorem{prop}[thm]{Proposition}
\newtheorem{lem}[thm]{Lemma}
\newtheorem{cor}[thm]{Corollary}

\theoremstyle{definition}
\newtheorem{defn}[thm]{Definition}

\newtheorem{ex}[thm]{Example}
\newtheorem{sit}[thm]{Situation}
\theoremstyle{remark}
\newtheorem{notn}[thm]{Notation}
\newtheorem{rem}[thm]{Remark}

\numberwithin{equation}{section}

%% file: twist-2.pdf_t
\begin{picture}(0,0)%
\includegraphics{twist-2.pdf}%
\end{picture}%
\setlength{\unitlength}{3947sp}%
\begingroup\makeatletter\ifx\SetFigFont\undefined%
\gdef\SetFigFont#1#2#3#4#5{%
  \reset@font\fontsize{#1}{#2pt}%
  \fontfamily{#3}\fontseries{#4}\fontshape{#5}%
  \selectfont}%
\fi\endgroup%
\begin{picture}(5324,2358)(-21,-3168)
\put(3834,-2856){\makebox(0,0)[lb]{\smash{{\SetFigFont{10}{12.0}{\rmdefault}{\mddefault}{\updefault}{\color[rgb]{0,0,0}$1$}%
}}}}
\put(424,-1780){\makebox(0,0)[lb]{\smash{{\SetFigFont{10}{12.0}{\rmdefault}{\mddefault}{\updefault}{\color[rgb]{0,0,0}$v$}%
}}}}
\put(1901,-1767){\makebox(0,0)[lb]{\smash{{\SetFigFont{10}{12.0}{\rmdefault}{\mddefault}{\updefault}{\color[rgb]{0,0,0}$1$}%
}}}}
\put(3304,-1763){\makebox(0,0)[lb]{\smash{{\SetFigFont{10}{12.0}{\rmdefault}{\mddefault}{\updefault}{\color[rgb]{0,0,0}$0$}%
}}}}
\put(2595,-3113){\makebox(0,0)[lb]{\smash{{\SetFigFont{10}{12.0}{\rmdefault}{\mddefault}{\updefault}{\color[rgb]{0,0,0}$0$}%
}}}}
\put(1353,-2856){\makebox(0,0)[lb]{\smash{{\SetFigFont{10}{12.0}{\rmdefault}{\mddefault}{\updefault}{\color[rgb]{0,0,0}$0$}%
}}}}
\put(4711,-1780){\makebox(0,0)[lb]{\smash{{\SetFigFont{10}{12.0}{\rmdefault}{\mddefault}{\updefault}{\color[rgb]{0,0,0}$v'$}%
}}}}
\put(1878,-2290){\makebox(0,0)[lb]{\smash{{\SetFigFont{10}{12.0}{\rmdefault}{\mddefault}{\updefault}{\color[rgb]{0,0,0}$0$}%
}}}}
\put(3317,-2279){\makebox(0,0)[lb]{\smash{{\SetFigFont{10}{12.0}{\rmdefault}{\mddefault}{\updefault}{\color[rgb]{0,0,0}$0$}%
}}}}
\put(2595,-933){\makebox(0,0)[lb]{\smash{{\SetFigFont{10}{12.0}{\rmdefault}{\mddefault}{\updefault}{\color[rgb]{0,0,0}$0$}%
}}}}
\end{picture}%

%% file: lls-non-ct-trop.bbl
\newcommand{\etalchar}[1]{$^{#1}$}
\newcommand{\noopsort}[1]{} \newcommand{\printfirst}[2]{#1}
  \newcommand{\singleletter}[1]{#1} \newcommand{\switchargs}[2]{#2#1}
\providecommand{\bysame}{\leavevmode\hbox to3em{\hrulefill}\thinspace}
\providecommand{\MR}{\relax\ifhmode\unskip\space\fi MR }
\providecommand{\MRhref}[2]{%
  \href{http://www.ams.org/mathscinet-getitem?mr=#1}{#2}
}
\providecommand{\href}[2]{#2}
\begin{thebibliography}{HLM{\etalchar{+}}08}

\bibitem[AB15]{a-b1}
Omid Amini and Matthew Baker, \emph{Linear series on metrized complexes of
  algebraic curves}, Mathematische Annalen \textbf{362} (2015), no.~1-2,
  55--106.

\bibitem[Bak08]{ba2}
Matthew Baker, \emph{Specialization of linear systems from curves to graphs},
  Algebra \& Number Theory \textbf{2} (2008), no.~6, 613--653.

\bibitem[BN07]{b-n1}
Matthew Baker and Serguei Norine, \emph{{R}iemann-{R}och and {A}bel-{J}acobi
  theory on a finite graph}, Advances in Mathematics \textbf{215} (2007),
  no.~2, 766--788.

\bibitem[BPR13]{b-p-r2}
Matthew Baker, Sam Payne, and Joseph Rabinoff, \emph{On the structure of
  non-archimedean analytic curves}, Tropical and Non-Achimedean Geometry,
  Contemporary Mathematics, vol. 605, American Mathematical Society, 2013,
  pp.~93--121.

\bibitem[CDPR12]{c-d-p-r}
Filip Cools, Jan Draisma, Sam Payne, and Elina Robeva, \emph{A tropical proof
  of the {B}rill-{N}oether theorem}, Advances in Mathematics \textbf{230}
  (2012), no.~2, 759--776.

\bibitem[CLM15]{c-l-m1}
Lucia Caporaso, Yoav Len, and Margarida Melo, \emph{Algebraic and combinatorial
  rank of divisors on finite graphs}, Journal de Math\'ematiques Pures et
  Appliqu\'ees \textbf{104} (2015), no.~2, 227--257.

\bibitem[EH86]{e-h1}
David Eisenbud and Joe Harris, \emph{Limit linear series: Basic theory},
  Inventiones Mathematicae \textbf{85} (1986), no.~2, 337--371.

\bibitem[GK08]{g-k2}
Andreas Gathmann and Michael Kerber, \emph{A {R}iemann-{R}och theorem in
  tropical geometry}, Mathematische Zeitschrift \textbf{259} (2008), no.~1,
  217--230.

\bibitem[He17]{he2}
Xiang He, \emph{Smoothing of limit linear series on curves and metrized
  complexes of pseudocompact type}, preprint, 2017.

\bibitem[HLM{\etalchar{+}}08]{h-l-m-p-p-w1}
Alexander Holroyd, Lionel Levine, Karola M\'ez\'aros, Yuval Peres, James Propp,
  and David Wilson, \emph{Chip-firing and rotor-routing on directed graphs}, In
  and Out of Equilibrium II, Progress in Probability, vol.~60, Birkhauser,
  2008, pp.~331--364.

\bibitem[Jen16]{je1}
David Jensen, \emph{The locus of {B}rill-{N}oether general graphs is not
  dense}, Portugaliae Mathematica \textbf{73} (2016), no.~3, 177--182.

\bibitem[KKW]{k-k-w1}
Sameer Kailasa, Vivian Kuperberg, and Nicholas Wawrykow, \emph{Chip-firing on
  trees of loops}, preprint.

\bibitem[LPP12]{l-p-p1}
Chang~Mou Lim, Sam Payne, and Natasha Potashnik, \emph{A note on
  {B}rill-{N}oether theory and rank-determining sets for metric graphs},
  International Mathematics Research Notices \textbf{2012} (2012), no.~23,
  5484--5504.

\bibitem[MZ08]{m-z1}
Grigory Mikhalkin and Ilia Zharkov, \emph{Tropical curves, their {J}acobians,
  and theta functions}, Curves and abelian varieties, Contemporary Mathematics,
  vol. 465, American Mathematical Society, 2008, pp.~203--230.

\bibitem[Oss]{os25}
Brian Osserman, \emph{Limit linear series for curves not of compact type},
  Journal f\"ur die reine und angewandte Mathematik (Crelle's journal), to
  appear.

\bibitem[Oss14]{os20}
\bysame, \emph{Limit linear series moduli stacks in higher rank}, preprint,
  2014.

\bibitem[Oss16]{os23}
\bysame, \emph{Dimension counts for limit linear series on curves not of
  compact type}, Mathematische Zeitschrift \textbf{284} (2016), no.~1-2,
  69--93.

\bibitem[Pfl]{pf2}
Nathan Pflueger, \emph{Special divisors on marked chains of cycles}, Journal of
  Combinatorial Theory, Series A, to appear.

\end{thebibliography}
